\documentclass[graybox]{svmult}

\usepackage{mathptmx}       
\usepackage{helvet}         
\usepackage{courier}        
\usepackage{type1cm}        
\usepackage{makeidx}         
\usepackage{graphicx}        
\usepackage{multicol}        
\usepackage[bottom]{footmisc}

\spdefaulttheorem{project}{Research Project}{\bf}{}
\newcommand{\resproject}[1]{\begin{svgraybox} \begin{project} #1 \end{project} \end{svgraybox}}

\spdefaulttheorem{cproblem}{Challenge Problem}{\bf}{}
\newcommand{\chproblem}[1]{\begin{cproblem} #1 \end{cproblem}}

\spnewtheorem*{prerequisites}{Suggested prerequisites}{\bf}{\small\em}

\makeindex             

\usepackage{amsmath}
\usepackage{amsfonts}
\usepackage{enumerate}
\usepackage{fancyvrb}

\begin{document}

\title*{What is the Best Way to Do Something? \\ A Discreet Tour of Discrete Optimization}
\titlerunning{A Discreet Tour of Discrete Optimization}
\author{Thiago Serra}
\institute{Thiago Serra \at University of Iowa, United States, \email{thiago-serra@uiowa.edu}
}
%
%
\maketitle


\abstract{
In mathematical optimization, 
we want to find the best possible solution for a decision-making problem. 
Curiously, 
these problems are harder to solve if they have discrete decisions. 
Imagine that you would like to buy chocolate: 
you can buy no chocolate or one chocolate bar, but typically you cannot buy just half of a bar. 
Now imagine that you could also buy many other items,  
and that you need to meet nutritional needs while minimizing the grocery bill. 
With more options and more demands, finding the best solution becomes trickier. 
But since many real-world settings benefit from mathematical optimization, 
such as scheduling trains and flights, planning truck deliveries, and making better investment decisions, 
these problems are widely studied in a branch of mathematics called Operations Research (OR). 
Sometimes we can simply write the mathematical model and find an optimal solution with OR software, but for larger problems we may need to develop new mathematical models and even write our own algorithms.  
We explore both cases  
with a simple and well-known problem (the traveling salesperson problem), some computer programming (in Python), and software that is free for academic use (Gurobi). 
}

\begin{prerequisites}
An introductory course in algorithms and computer programming, preferably but not necessarily in Python, and basic understanding of calculus and linear algebra. 
Some familiarity with statistical visualization tools may be helpful for analyzing experiments. 
\end{prerequisites}

\section*{Preface}

I am writing this text as a friendly introduction to mathematical optimization for undergraduate students. 
The purpose is to provide the tools for students to explore research projects in this area.
I opted to write it as close to how I would have lectured about it, 
which does not read as a typical academic paper. 
In fact, I start by assuming very little background from the reader. 
But as you read along, 
you will see more interesting applications and also academic references for further reading toward the last sections. 
For helping we understand how to talk about optimization, 
I would like to acknowledge every undergrate student who did research with me at Bucknell University, 
and also every student who took my optimization lectures at Carnegie Mellon University, Bucknell University, and the University of Iowa. 
Thank you!

I also would like to thank Allison Lewis, Aaron Wootton, Alex John Quijano, and Richard Forrester for their feedback on a preliminary version of this text.

All the code and data used is available at: \url{https://github.com/thserra/discreet}

\section{Introduction}\label{sec:intro} 

In mathematical optimization, 
our goal is to find the best possible solution for a decision-making problem. 
We define as a solution the values that we choose for variables representing the decisions, 
which are the \emph{decision variables}. 
For example, imagine that we want to decide how many apples and oranges to buy. 
Let us call these decision variables $a$ (for apples) and $o$ (for oranges), and then write a mathematical model with them. 
We start with a formula for what we would like to optimize, 
which is the \emph{objective function}. 
For example, imagine that each apple costs \$0.49 and that each orange costs \$0.50. 
We express cost minimization as follows:
\[
\min ~~~ 0.49 a + 0.50 o
\]

But why are we buying fruit? 
If our goal is to minimize cost, then the best thing to do is to buy nothing. 
Better yet, in case we already have some apples and oranges, 
we could even make some money by selling them. 
For example, 
we could say that $a=1$ means that we are buying one apple, 
whereas $o=-2$ means that we are selling two oranges. 
But if we do not have apples or oranges of our own, or if we do not want to sell the ones that we have, 
then we should make that clear in our model. 
When we want to limit the options available, 
we define a \emph{constraint}. We can say: 
\[
a \geq 0 \]
\[
o \geq 0
\]
Or, to make things simpler, we may use one-liners like the following one:
\[
a, o \geq 0
\]
In fact, this is a special type of constraint because it clearly defines the minimum value of the decision variables. We call it a \emph{lower bound}. 
If we were defining the maximum value of the decision variables, 
we would have an \emph{upper bound}. 

At this point, we have said that fruit costs money and that we are not selling any fruit. 
Hence,  
the optimal solution would be to buy nothing because we would spend nothing. 
Unless, of course, we have a reason to buy fruit. 
In fact, 
if we did not, 
it would have been a waste of time to be talking about it! 
So imagine that you would like to buy fruit to supply at least 5\% of your daily need of potassium\footnote{To make the problem more realistic, daily values of nutritional needs were obtained from the FDA website~\cite{fda2024nutrition} and fruit nutrients were obtained from the USDA website~\cite{usda2024apple,usda2024bananas,usda2024oranges,usda2024pears}.}. 
With each apple having 214 mg of potassium, each orange having 232 mg of potassium, and 235 mg being 5\% of the daily need of potassium, we use the following constraint:
\[
214 a + 232 o \geq 235
\]

Now that our problem is finally nontrivial, what would be an optimal solution? 
This would be a good time to pause and think about it. 
After all, thinking about what we expect an optimal solution to be may help us figure out if something is missing. 
That is what we did with the objective function, which led to the lower bounds,  
and now to a constraint on potassium. 
The same strategy works well with computer programming: 
write a small piece that does some work, understand if it really does what you want it to do, and then add and test other small pieces until you are done.  

So I would recommend that you stop reading right now if you have not thought about how to solve the problem that we have formulated so far. 
Seriously, if you have not done that yet, stop and think. I am tired of adding fillers already!

It is tempting to look one paragraph ahead, isn't it?

Perhaps you thought that you could see the answer two paragraphs ahead?

A single apple or a single orange would not have enough potassium, 
but we would have enough potassium with two apples, two oranges, or a combination of one apple and one orange. 
Apples are cheaper, so we could choose $a=2$ and $o=0$. 

However, the solution of our model would be $a = 0$ and $o \approx 1.013$, 
so you would not have gotten away with the answer even by looking three paragraphs ahead!

But why did we go from whole apples to sliced oranges? 
First, we did not specify that the decision variables can only take integer values. 
We usually assume that decision variables are continuous unless stated otherwise, which can be confusing.  
Hence, even if that is the case, it is a good practice to make it explicit. 
For example, we can use $a, o \in \mathbb{R}$ if the variables should be continuous, or 
$a, o \in \mathbb{R}_+$ if we also want to exclude negative values. 
For nonnegative integer variables, we use:
\[
a, o \in \mathbb{Z}_+
\]
Now we are actually defining the \emph{domain} of the decision variables, 
which is the set of all values that they can take: any integer value greater than or equal to zero.

Second, if we do not have to buy whole fruit, 
then the problem becomes finding the cheapest source of potassium per dollar, regardless of how much we need. 
With apples, you get $\frac{214}{0.49} \approx 398$ mg of potassium per dollar. 
With oranges, you get $\frac{232}{0.50} \approx 458$ mg of potassium per dollar. 
Hence, buying $\frac{235}{232} \approx 1.013$ oranges (for $\approx \$0.51$) would be the cheaper than buying $\frac{235}{214} \approx 1.010$ apples (for $\approx \$0.54$).

This example is important because it shows that we cannot get away with rounding continuous variables to the nearest integer values. 
That is a tempting route if you consider that solving a problem with continuous variables tends to be easier, 
especially if the objective function and the constraints are all linear. 
Curiously, pretending that the decision variables are all continuous is a useful step in solving discrete optimization problems. 
This is where we will go in the next two sections. 

\begin{exercise}
In the example above, 
the optimal solution with integer variables is to buy two apples for $\$0.98$. 
However, the optimal solution with continuous variables is to buy $\approx 1.013$ oranges for $\approx \$0.51$, 
which is much cheaper. 
We call the value of the objective function in an optimal solution as the \emph{optimal value}. 
If the only difference between two minimization models $A$ and $B$ is that the decision variables are continuous in $A$ and integer in $B$, 
is the optimal value for $A$ always smaller?
\end{exercise}

\begin{exercise}
In the example above, 
the optimal solutions of the models with continuous and with integer decision variables are so different on purpose.  
That is how instructors create a ``teachable moment'' and help students understand the nuance of what we are discussing. 
Now it is your turn to change the model on purpose:  
\begin{enumerate}[(a)]
\item If you can only change the potassium need, say from $5\%$ of the daily need to some other percentage, how can you make the solution of both models identical?
Are there other choices for that percentage that would have the same effect?
\item Now let us consider only the model with integer variables, but let us ignore the changes that you may have made when answering item (a).
By only changing the price of oranges now, and by changing it by the least amount possible to make two oranges an optimal solution, what would be the new price?
\item How many solutions would be optimal with such a change to the orange price? 
\end{enumerate}
\end{exercise}

\begin{svgraybox}
This is a good point to mention that optimizers are often careful with what they say. 
For example, they may avoid talking about \emph{the} optimal solution 
because some problems have more than one optimal solution. 
Sometimes we do know that a problem has only one optimal solution, 
which we call the \emph{optimum}; 
and other times we know or assume that a problem has multiple optimal solutions, 
which we call the \emph{optima}, 
in which case we are indifferent among them.  
But we never, ever, say that one solution is ``more optimal'' than another!
\end{svgraybox}

\begin{exercise}
Continuing with the same example, 
what if we already had an orange? 
If we were to keep the orange, then we would cover some of the potassium need. 
If we were to sell the orange, then we would offset part of the cost of buying apples. 
\begin{enumerate}[(a)]
\item Assume first that we can buy and sell fruit at the same price. If we keep $a$ and $o$ as decision variables for how many apples and oranges we buy (with a negative value implying that we sell), how should we change the objective function, constraint, and bounds to model the case in which we already have an orange?
\item Now suppose that we only get paid $\$0.48$ for our orange, even though we have to pay $\$0.50$ to buy another one. This is what happens when you trade currency, say US dollars for Brazilian reais: the shop buys for less than they sell, since maintaining a shop costs money and they also need to profit for a living. Digression apart, how can you model having an orange to sell along with buying apples and oranges, such that we spend the least to meet the potassium need?
\end{enumerate}
\end{exercise}

\section{Models}\label{sec:model}

So far we have seen only a small model in parts. 
Now we will see a whole model, with a few more decision variables and constraints. 
First, we now choose among buying apples, oranges, bananas, and pears. 
That is not a lot of fruit yet, 
but this is a good time to stop using separate variables for each fruit. 
Otherwise, we would run out of letters in the alphabet at some point. 
Hence, we use $x_i$ for how much we buy of fruit $i$, 
where $i=1$ is for apples, $i=2$ for oranges, $i=3$ for bananas, and $i=4$ for pears. 
Second, besides buying enough fruit to meet 5\% of the daily need of potassium (235 mg), 
we also want to meet 5\% of the daily need of calcium (65 mg) and 20\% of the daily need of fiber (5.6 g). 
Finally, we price bananas at \$0.51 and pears at \$0.52:
\[
\begin{aligned}
\min ~~& 0.49 x_1 + 0.50 x_2 + 0.51 x_3 + 0.52 x_4 \\
\text{s.t.} ~~&  214 x_1 + 232 x_2 + 375 x_3 + 155 x_4 \geq 235 \\
&  12 x_1 + 60.2 x_2 + 5.75 x_3 + 14.2 x_4 \geq 65 \\
&  4.8 x_1 + 2.8 x_2 + 5.31 x_3 + 5.52 x_4 \geq 5.6 \\
& \textbf{x} \in \mathbb{Z}_+^4
\end{aligned}
\]
The first line of a mathematical optimization model is typically the objective function, preceded by either $\min$ (for \emph{minimize}) or $\max$ (for \emph{maximize}). 
The second line is typically the first one with a constraint and the one in which s.t. (for \emph{subject to}) is written. 
After the constraints connecting multiple decision variables, 
the last lines usually have the bounds or the domain of each decision variable. 
Above, we have $\textbf{x}$ as a vector of 4 variables with nonnegative and integer values. 

Now we have a problem that is definitely not so straightforward to solve, 
but we will wait until Section~\ref{sec:bb} to talk about how to solve it. 
Before that, 
we will go in more detail about different types of models, why they matter, and how they relate to one another. 
This will be useful to understand how we solve optimization problems.

\subsection{The General Form}

Our new model is another example of an \emph{Integer Linear Programming}~(\textbf{ILP}) model. 
We say that this is an \emph{integer} model because all the decision variables need to have integer values. 
We say that this is a \emph{linear} model because the contribution of each decision variable to the objective function and to the constraints is linearly proportional to its value. 
For example, the decision variables should not be multiplied by one another (as in $x_1 x_2$) or used as the argument of a nonlinear function (as in $\cos(x_1)$) 
in a constraint or the objective function. 
We generally prefer to formulate models that are linear 
because there are faster algorithms to solve them. 
We say that this is a \emph{programming} model because we use it to make plans. 
Although we often use computer programming to work with models, 
this has nothing to do with computer programming. 
In fact, some practitioners prefer calling it a  
\emph{prescriptive model}\footnote{In business schools, the application of mathematical models and computers to solve applied problems based on data is often called \emph{business analytics}, 
which is subdividided into \emph{descriptive analytics} for the application of statistical methods to describe historic data, \emph{predictive analytics} for the application of machine learning and statistical models to produce forecasts, and \emph{prescriptive analytics} for the application of mathematical optimization to make decisions.  
We have avoided using the term \emph{planning} because it is used by a traditional branch of artificial intelligence focused on agent decision-making tasks, such as making a robot navigate a room from one place to another.} for clarity. 

We can write an ILP model with $n$ decision variables in $\textbf{x}$ and $m$ constraints as 
\[
z = \min \{ \textbf{c}^T \textbf{x} ~ |  ~ A \textbf{x} \geq \textbf{b},  ~ \textbf{l} \leq \textbf{x} \leq \textbf{u}, ~ \textbf{x} \in \mathbb{Z}^n \}, 
\]
where $z$ is a variable commonly used for the value of optimal solutions. 
This formulation also includes vector 
$\textbf{c} \in \mathbb{R}^n$ for costs, matrix $A \in \mathbb{R}^{m \times n}$ for the left-hand side (LHS) coefficients of the constraints, vector $\textbf{b} \in \mathbb{R}^m$ for the right-hand side (RHS) of the constraints, and vectors $\textbf{l}, \textbf{u} \in \mathbb{R}^n$ for the lower and upper bounds on the values of the decision variables. 
For example, if we were to write down the first constraint explicitly, we would have 
\[
a_{1 1} x_1 + a_{1 2} x_2 + \ldots + a_{1 n} x_n \geq b_1
\]
or 
\[
\sum_{i=1}^n a_{i j} x_j \geq b_1.
\]

One special type of integer decision variable that we may see very often is the \emph{binary} decision variable, which only takes the values 0 and 1.
For example, we may use 
\[
x \in \{0, 1\}
\]
to denote that $x$ is a binary variable, whereas we may use 
\[
\textbf{x} \in \{0, 1\}^3
\]
to denote that $\textbf{x}$ is a vector of three binary variables: $x_1$, $x_2$, and $x_3$. 

\subsection{Models vs. Problems}

So far we have talked a little about optimization problems and a lot about optimization models. 
By a model, 
we consider a mathematical formulation for an optimization problem. 
Some prefer to talk about an ILP formulation instead of an ILP model, but they mean the same thing. 
However, models and problems are not the same. 
On the one hand, 
we may come up with different models to solve the same problem, 
as we will talk more about in Section~\ref{sec:tricks}. 
In fact, 
there is a common mistake of confusing problem description with formulation. 
It is important to understand that defining the problem is not the same thing as expressing it mathematically with decision variables, objective function, and constraints: that comes later. 
On the other hand, 
changing a model may imply that we are solving a different problem, 
as we will see next.

\subsection{Other Types of Models}

In discrete optimization, there are some important variations of the ILP model. 
The first one is the \emph{Linear Programming} (\textbf{LP}) model, 
in which all decision variables are continuous and that last segment denoting the domain of the variables becomes 
\[
\textbf{x} \in \mathbb{R}^{n}.
\]
Isn't it strange that LP, a type model with no discrete decision variables, is important for discrete optimization? 
The main reason is that it tends to be computationally easier to solve an LP model than it is to solve an ILP model. 

If we start with an ILP model but change the domain of the decision variables to be continuous, 
from $\mathbf{x} \in \mathbb{Z}_+$ to $\mathbf{x} \in \mathbb{R}_+$, 
then we have an LP model that is a \emph{relaxation} of the ILP model. 
We say that the LP model is a relaxation because 
all the solutions of the ILP model are \emph{feasible}\footnote{For historic reasons, optimizers say \emph{feasible solution} when they mean that the solution is valid, in the sense that it satisfies all the constraints of the model---including bounds and domains.} for the LP model, but not all solutions of the LP model are feasible for the ILP model. 
For example, 
in that first model in which we could only buy apples and oranges, 
the solution $(a, o) = (2, 0)$ is feasible for both models,  
whereas the solution $(a, o) = (0, 1.013)$ is only feasible for the LP relaxation---since the ILP model does not admit slicing fruit. 
Hence, it is as if we were solving a slightly different problem, 
one that might be a little easier to solve because we have more options. 
Sometimes we are lucky enough that an optimal solution that we find for the LP model also works for the ILP---and that is just one way in which we can benefit from LP relaxations. 
We will see this in more detail in the next section. 

The other important variation of the ILP model is the \textbf{Mixed--Integer Linear Programming} (\textbf{MILP}) model, 
in which some variables are integer and other variables are continuous, 
and---\emph{without loss of generality}\footnote{Without loss of generality, or w.l.o.g., is a curious mathematical term. Basically, we are acknowledging that what happens next looks like a special case but, if you really think about it, it is not. For example, if we say that $p$ out of $n$ decision variables are integer, then they may not be necessarily the first $p$ decision variables. However, it is convenient for us to describe an MILP model this way, and so we use w.l.o.g to imply that changing the order of the variables would not affect the essence of what we are discussing.}---the domain becomes 
\[
\textbf{x} \in \mathbb{Z}^p \times \mathbb{R}^{n-p}
\]
with the first $p$ decision variables being integer, $x_1$ to $x_p$, and the remaining $n-p$ decision variables being continuous, $x_{p+1}$ to $x_{n}$. 
When we talk about a relaxation of an MILP model, 
we are talking about relaxing the first $p$ integer variables to be continuous. 
We will not discuss many MILP models in this text, 
but you may come across more MILP models if you continue studying discrete optimization.

\begin{exercise}\label{ex:knapsack}
Let us assume that Table~\ref{tab:1} describes the author's happiness from consuming different chocolate bars. 
Each of those bars has a different price. 
This is an example of the \emph{knapsack problem}, 
in which the goal is to maximize the benefit from selecting among a collection of items 
while not exceeding their weighted sum.
\begin{enumerate}[(a)]
\item Formulate an ILP model in which $x_i$ is how many bars of chocolate of type $i$ we buy, 
subject to a budget of \$10.00, 
and to a maximum of one bar of each type. 
\item What is an optimal solution of your ILP model?
\item What is an optimal solution of the LP relaxation of your model?
\item If we can buy as many bars of chocolate of each type as we want, 
how would that change your answers to items (a), (b), and (c)?
\end{enumerate}
\end{exercise}

\begin{table}
\centering
\caption{Data for the knapsack problem to be modeled in Exercise~\ref{ex:knapsack}.}
\label{tab:1}       
%
%
\begin{tabular}{cccccc}
\textbf{Chocolate} & Purity & Lindt & Dove & Reese's & Hershey's \\
\hline\noalign{\smallskip}
\textbf{Number} & 1 & 2 & 3 & 4 & 5 \\
\textbf{Price} & \$4.50 & \$4.00 & \$3.00 & \$3.00 & \$2.00 \\
\textbf{Happiness} & 10 & 8 & 7 & 5 & 1  
\end{tabular}
\end{table}

\chproblem{Is there a redundant constraint in the model described in Section~\ref{sec:model}? A constraint is not redundant if, and only if, there are solutions that are infeasible only because that constraint exists.}

\chproblem{\label{ex:facility}
A supermarket chain wants to rent storage facilities from where their stores will be supplied. 
Every store will be supplied by one storage facility, but a storage facility may supply multiple stores. 
There are $m$ possible storage facilities to be rented, numbered from $i = 1$ to $i = m$. 
If store facility $i$ is rented, the monthly rent would be $r_i$ and we would have $s_i$ of space in that storage facility.
The supermarket chain has $m$ stores, numbered from $j = 1$ to $j = n$. 
The total space required by goods that will be shipped to store $j$ is $d_j$. 
If storage facility $i$ supplies store $j$, there will be a monthly cost of $t_{i j}$.\footnote{Note how we consistently use index $i$ when referring to a facility and index $j$ when referring to a store. That is a good practice to avoid inconsistencies in the model and bugs in the code. For example, we expect index $i$ to range from 1 to $m$ and index $j$ to range from 1 to $n$, so it is less likely that we would write a model or code in which we mistakenly assume to have $n$ (instead of $m$) facilities or $m$ (instead of $n$) stores, which would be a problem --- especially if $m \neq n$.} 
This is an example of the \emph{facility location problem}. 

To model this problem, 
let us assume that we use a binary decision variable $x_{i j} \in \{0, 1\}$ denoting if storage facility $i$ supplies store $j$, 
and therefore we have $x \in \{0, 1\}^{m \times n}$, 
so that we can represent the monthly transportation cost as $\sum\limits_{i=1}^m \sum\limits_{j=1}^n t_{i j} x_{i j}$. 
Moreover, 
let us assume that we use a binary decision variable $y_i$ denoting if store facility $j$ is rented, 
so that we can represent the monthly rent as $\sum\limits_{i=1}^m r_i y_i$.
\begin{enumerate}[(a)]
\item Write a complete model for this problem, so that the supermarket chain minimizes their monthly costs while satisfying the following conditions:
\begin{itemize}
\item Every store should be assigned to one, and only one, storage facility.
\item The sum of the space needed by all the stores assigned to a given storage facility should not exceed the space available in that storage facility.
\item If store $j$ is assigned to storage facility $i$ ($x_{i j} = 1$), then storage facility $i$ should be rented ($y_i = 1$).\footnote{This last condition is a tricky one. If you get stuck, keep reading and we will talk about it later.}
\end{itemize}
\item If the supermarket did not have to pay rent for the storage facilities and every storage facility had enough space to supply all the stores, 
what would be an optimal solution?
\end{enumerate}
}

\begin{svgraybox}
Some details make the difference between an optimization problem that is easy to solve and an optimization problem that is hard to solve. 
The part (b) of Challenge Problem~\ref{ex:facility} presents a special case in which the facility location problem becomes easy to solve. 
We will go back to the distinction between easy and hard problems at the end of the next section.
\end{svgraybox}

\section{The Branch-and-Bound Algorithm}\label{sec:bb}

In life, 
one way to deal with a problem is to pretend that the problem does not exist. 
Some people call this ``hide and hope''---while they hide, they hope that the problem gets sorted out on its own! 
Optimization problems do not get solved on their own, 
but sometimes we can benefit from ignoring some of their complicating parts. 

For example, we know that optimization problems with continuous decision variables are not as difficult to solve.\footnote{Because we focus on discrete optimization and this is supposed to be a brief introduction, 
we do not discuss algorithms to find optimal solutions of LP models. For a friendly introduction, see \cite{co_book}. 
For a shorter introduction like this text, see \cite{kotas}.}
So we may ignore for a moment that we have integer decision variables, 
which is why we talked about the LP relaxation previously. 
If we are lucky enough, 
we might find a solution that happens to have integer values only. 
In that case, at least the integrality problem got sorted out on its own. 

Branch-and-bound is what we typically do when that does not work. 
More precisely, branch-and-bound is what most OR software eventually do when a problem is more difficult to solve. 
It is something that is interesting at first, but that gets very tedious after a while. 
However, if we understand how it works and what makes it necessary, 
then we can design better models and avoid branching too much. 

To describe how branch-and-bound works, 
we will use it to solve the four fruit model from Section~\ref{sec:model}. 
If we solve the LP relaxation of that model, 
we find an optimal solution $(x_1, x_2, x_3, x_4) \approx (0.00, 0.95, 0.00, 0.53)$ with the optimal value $z \approx 0.75$. 
That is not a feasible solution for the ILP model. 
However, it shows that no solution has a better value than $\$0.75$. 
This is a \emph{lower bound on the optimal value}. 

If all the decision variables should be integer, 
then we know that there should be no decision variable $x_i$ such that $0 < x_i < 1$ in a feasible solution. 
Hence, 
if a solution exists, 
then either $x_i \leq 0$ or $x_i \geq 1$ for any $i$.  
So we could break a problem for which we cannot find a solution with our model 
into two smaller problems, 
one in which we add the constraint $x_i \leq 0$ and another one in which we add the constraint $x_i \geq 1$, 
knowing that any optimal solution of the original problem will also be an optimal solution of one of those models for the smaller problems. 

\begin{figure}[h]
\sidecaption
\includegraphics[scale=.23]{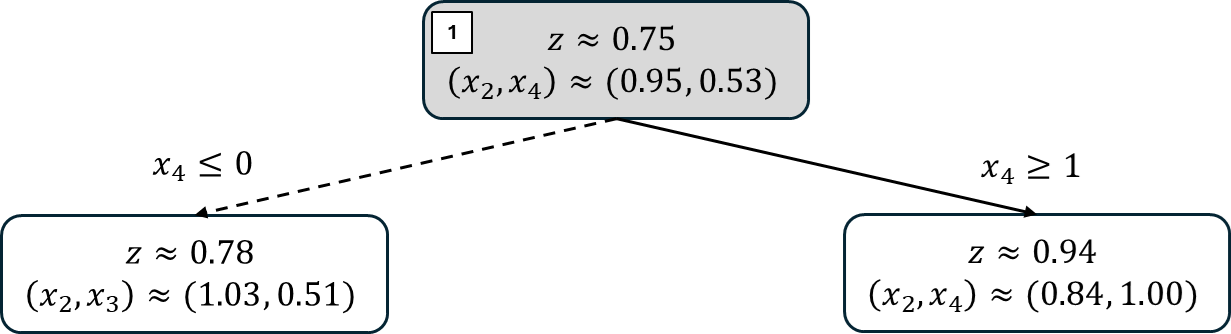}
%
%
\caption{Exploration of the root node (labeled 1) of the branch-and-bound tree by constraining decision variable $x_4$, which increases the lower bound of the optimal value from 0.75 to 0.78.}
\label{fig:BnB_01}       
\end{figure}

As a first step, let us pick variable $x_4$, 
which has the furthest away value from integer values in the solution of the LP relaxation, 
and then solve the subproblems shown in Figure~\ref{fig:BnB_01}, 
corresponding to buying no pear ($x_4 \leq 0$) and buying at least one pear ($x_4 \geq 1$).\footnote{Choosing the most fractional variable is a common tactic, 
but we could also have chosen to define subproblems corresponding to buying no orange ($x_2 \leq 0$) and buying at least one orange ($x_2 \geq 1$).}

Now that we have turned one problem into two, 
we have started building a \emph{branch-and-bound tree}, 
as illustrated in Figure~\ref{fig:BnB_01}. 
This tree has a root node representing the model for the original problem, 
which is gray and has a number 1 in it. 
Inside this node we have the optimal value of the LP relaxation 
and the value of the variables that are not zero in the LP solution, 
which we have seen before.
This node also has two children, 
representing models for subproblems in which we restrict the number of pears to be bought, 
as mentioned earlier.  

By solving the models of the children nodes, 
we keep finding noninteger values and still have no feasible solution. 
However, 
we find out that there is no feasible solution better than $\$0.78$, 
which is the optimal value of the LP relaxation of the left child, 
since the optimal value of the LP relaxation of the right child is greater than that.  
Hence, 
we can update the lower bound to $0.78$, 
which is more accurate. 

\begin{figure}[h]
\centering
\includegraphics[scale=.235]{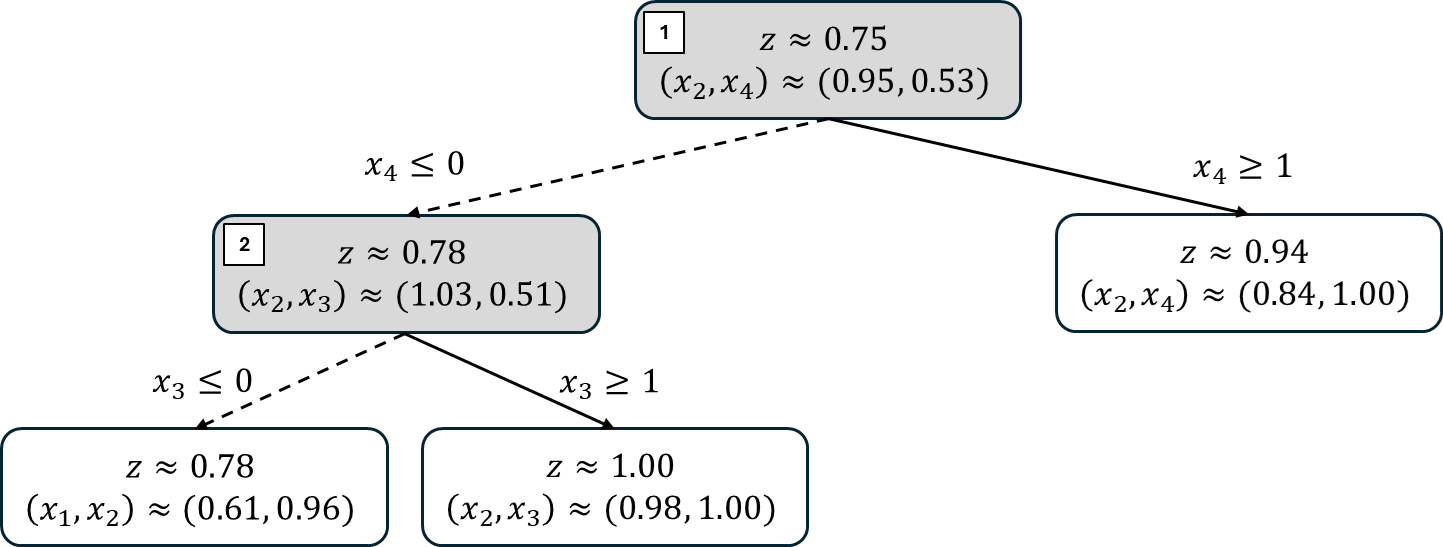}
\caption{Exploration of node 2 (left branch of the root node) of the branch-and-bound tree by constraining decision variable $x_3$, which does not change the lower bound significantly.}
\label{fig:BnB_02}       
\end{figure}

Having tried to crack the problem with one stab, 
here we stand, 
poor fools, 
barely wiser than before!\footnote{There is a famous quote by Johann Wolfgang von Goethe that is translated as 
``I HAVE, alas! Philosophy,
Medicine, Jurisprudence too,
And to my cost Theology,
With ardent labour, studied through.
And here I stand, with all my lore,
Poor fool, no wiser than before.''~\cite{faust}} 
So what do we do next? 
We keep breaking it apart, 
as in Figure~\ref{fig:BnB_02}.

We can choose either child of node 1 to proceed. 
We pick the left child, 
which now has the number 2 in Figure~\ref{fig:BnB_02}, 
and branch with variable $x_3$. 
That leads to two new subproblems in the children of node 2: 
the left one in which we do not buy pears or bananas ($x_4 \leq 0$ and $x_3 \leq 0$), 
and the right one in which we do not buy pears but we do buy bananas ($x_4 \leq 0$ and $x_3 \geq 1$). 
The lower bound is now the smallest optimal value among the nodes other than 1 and 2. Sadly, it grows less than a cent.

\begin{figure}[h]
\centering
\includegraphics[scale=.235]{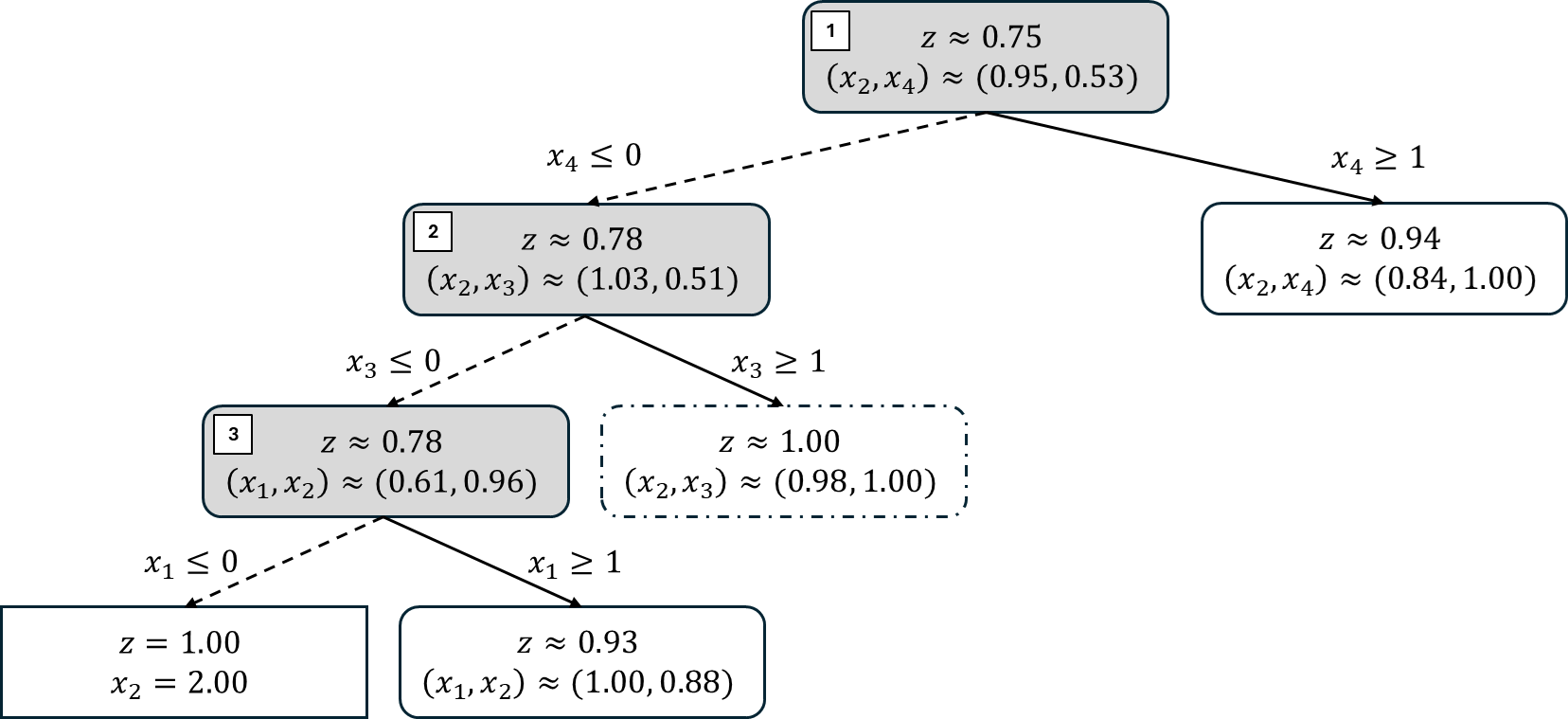}
\caption{Exploration of node 3 (left branch of node 2) of the branch-and-bound tree by constraining decision variable $x_1$, which produces an incumbent solution of value 1.00 (left branch), fathoms an unexplored node (right branch of node 2), and increases the lower bound from 0.78 to 0.93.}
\label{fig:BnB_03}       
\end{figure}

The nodes in gray (1 and 2) are \emph{explored}, 
meaning that we have already branched on them. 
The other nodes are \emph{unexplored}. 
The unexplored nodes are the only ones that matter going forward, 
since we will find from them and their children what is an optimal solution for the problem. 
We continue by exploring the leftmost child of node 2 in Figure~\ref{fig:BnB_03} by branching with $x_1$. 

This time we are luckier.  
From the third explored node in Figure~\ref{fig:BnB_03}, 
the left child (a rectangle) has a feasible solution of buying two oranges for $\$1.00$. 
Hence, 
we now have an \emph{upper bound on the optimal value}. 
The right child of node 3 happens to produce a new lower bound of $\$0.93$, 
since it becomes the smallest optimal value of the LP relaxation among all unexplored nodes. 
Moreover, 
we can now discard the right child of node 2, 
since we would not be able to find a better solution from it than the one that we already have. 
We say that this node has been \emph{fathomed}. 

The bounding step, 
through which we reduce the gap between the upper and the lower bounds, 
is what prevents us from doing an \emph{almost}\footnote{We will return to this point in Exercise~\ref{ex:enumeration}. \label{ft:enumeration}} complete enumeration of feasible solutions with branch-and-bound. 
Still, we have $\$0.07$ left between the best solution found and the theoretical best solution at this point, 
so we must proceed!

\begin{figure}[h]
\centering
\includegraphics[scale=.235]{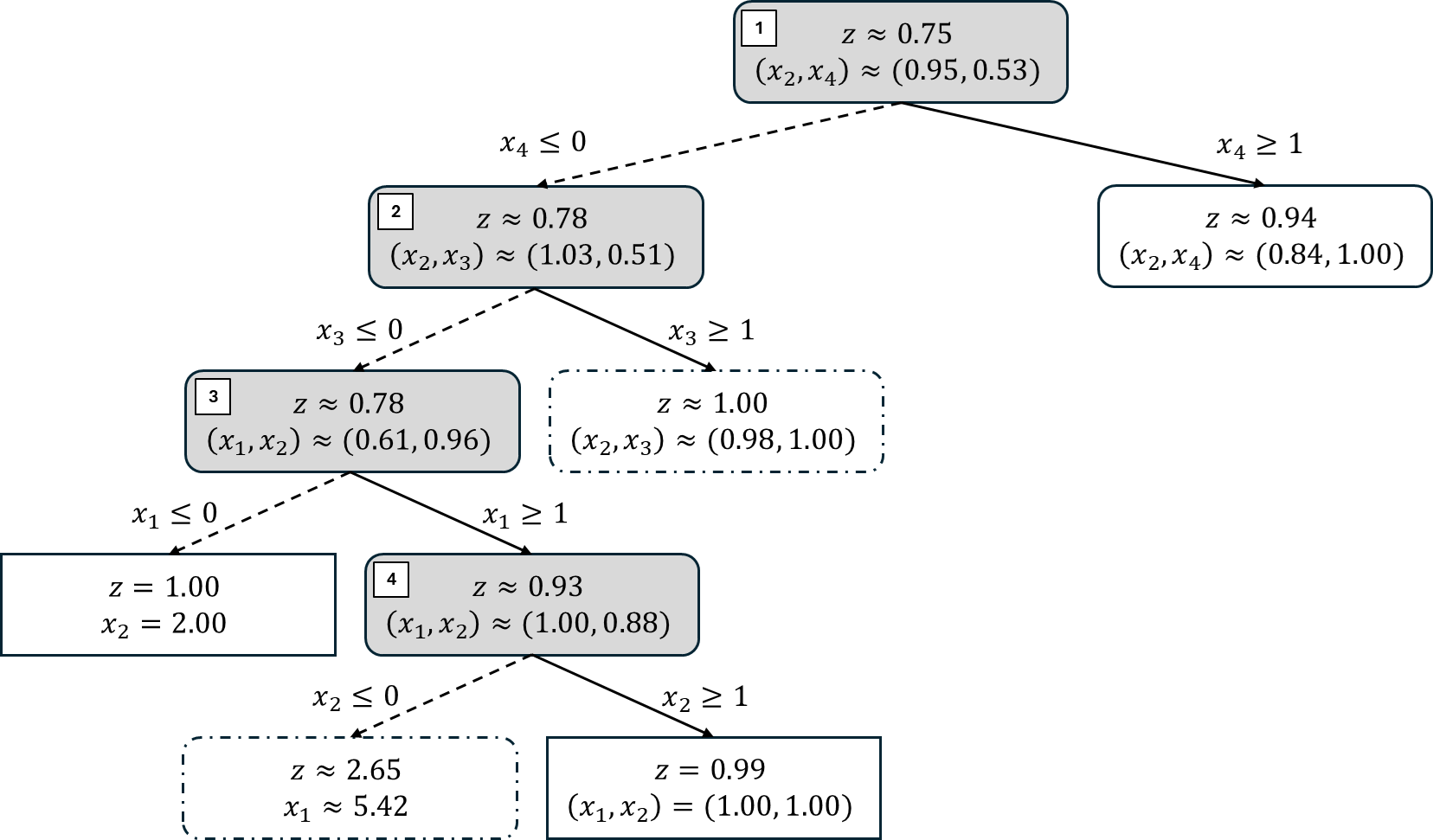}
\caption{Exploration of node 4 (right branch of node 3) of the branch-and-bound tree by constraining decision variable $x_2$, which produces an incumbent solution of value 0.99 (right branch), fathoms an unexplored node (its left branch), and increases the lower bound from 0.93 to 0.94.}
\label{fig:BnB_04}       
\end{figure}

We branch on the right child of node 3 in Figure~\ref{fig:BnB_04}, 
through which we obtain a left node that is automatically fathomed due to $z$ 
and a right node producing a better feasible solution of buying one apple and one orange for $\$0.99$, 
hence lowering the upper bound.  
As we have explored the node with a lower bound of $\$0.93$, 
we raise the lower bound to the value of the only unexplored node at this point: $\$0.94$. 

\begin{figure}[h]
\centering
\includegraphics[scale=.235]{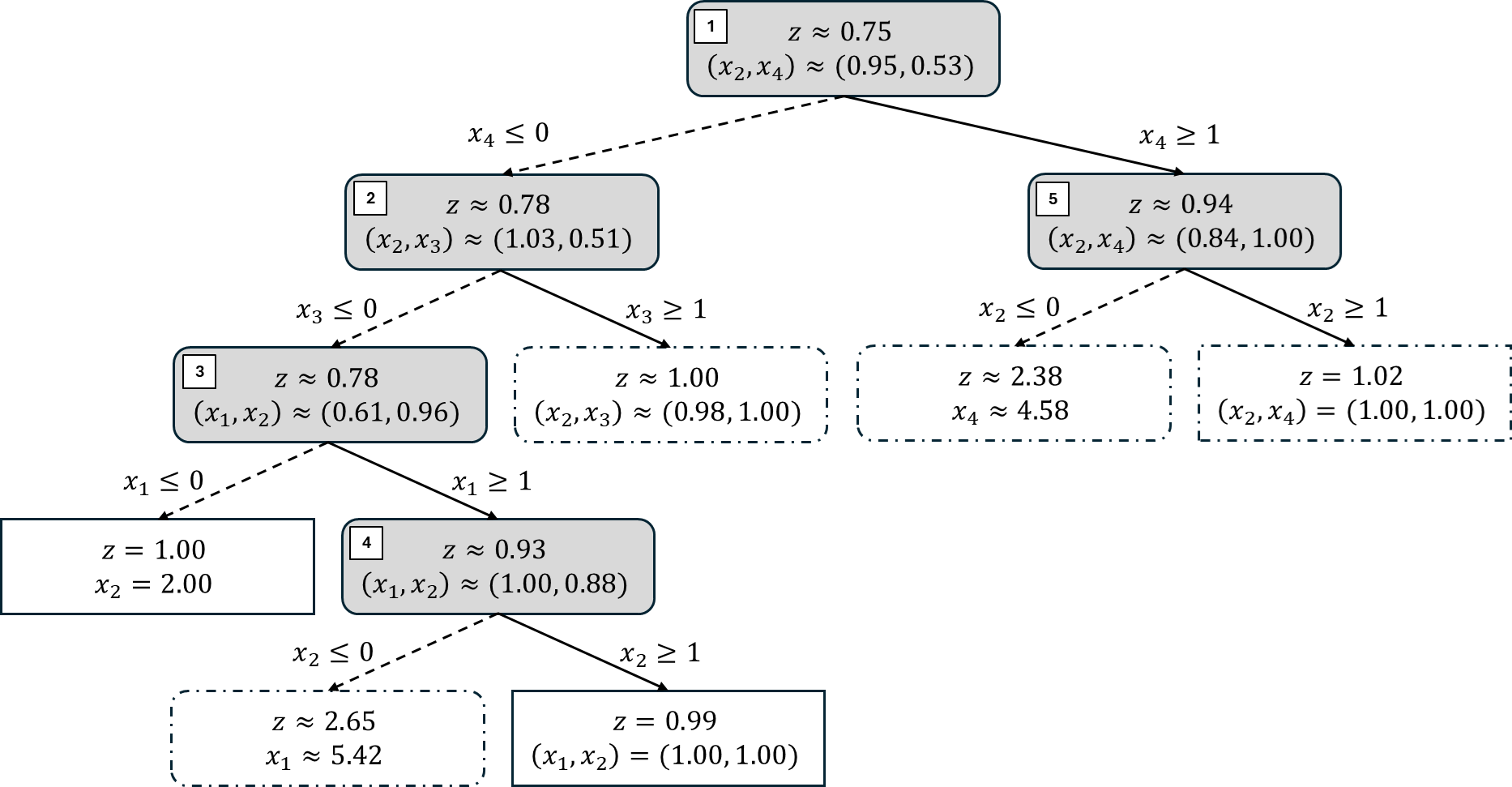}
\caption{Exploration of node 5 (right branch of node 1) of the branch-and-bound tree by constraining decision variable $x_2$, which produces two fathomed nodes (its branches), increases the lower bound from 0.94 to 0.99, and concludes that the incumbent solution of value 0.99 is optimal.}
\label{fig:BnB_05}       
\end{figure}

Now we have only one unexplored node (the right child of node 1) with only one noninteger variable to branch with ($x_2$), 
which when explored becomes node 5 in Figure~\ref{fig:BnB_05}. 
We obtain another feasible solution from it, 
but both child nodes are fathomed due their values. 
Since we do not have unexplored nodes left, 
we conclude that the optimal value is $\$0.99$, 
with $(x_1, x_2, x_3, x_4) = (1, 1, 0, 0)$ as an optimal solution. 

\begin{svgraybox}
\textbf{Where is the line between easy and hard optimization problems?} \\
We say that a problem is ``easy'' if we know of an algorithm that can solve it in time polynomial on the size of the problem. 
In our case, the size would involve the number of constraints $m$, the number of decision variables $n$, and the number of integer decision variables $p$. 
That is the case of solving LP models, for which reason we say that they are in the P class. 
In the case of ILP, we simply do not know if such an algorithm exists, even for just finding a feasible solution instead of an optimal one.  
In other words, the algorithms that we know for solving ILP models in general may take a runtime that eventually exceeds any polynomial time function on the size of the problem. 
That is why we call them ``hard''. 
In the case of solving ILP models, 
finding a feasible solution belongs to the NP-complete class 
and finding an optimal solution belongs to the NP-hard class.  
Those classes have many other problems in them.  
More interestingly, 
if we were able to find an algorithm that solves any of these problems in polynomial time, 
then it would be possible to solve all of them in polynomial time. 
That is the famous ``P = NP'' question: 
finding such an algorithm would imply that all of these problems are in P, 
and proving that no such algorithm exists would imply that NP-complete and NP-hard problems are indeed harder to solve.  
This is an open problem yet.
\\
\textbf{Does that matter?} \\
Not so much. 
Even though most optimizers believe that no such algorithm exists, 
they will not shy away from an NP-hard problem 
because there are many practical applications for those. 
This is what most research in mathematical optimization is about: 
we might not be able to solve the models fast enough when they exceed a certain size, 
but what really matters is if we can solve them for the sizes that we find in practice. 
That is the real impact.
\end{svgraybox}

\begin{exercise}
A lot has gone unmentioned in our example of how branch-and-bound works. 
We ask that you think and provide a short answer to the following questions: 
\begin{enumerate}[(a)]
\item Buying two apples as in an optimal solution of our first model in Section~\ref{sec:intro} would have been cheaper than $\$0.99$. Why can't we do that now?
\item By only branching on the values of $x_4$, $x_3$, and then $x_1$, 
why is it that the left child of node 3 produced a solution with an integer value for $x_2$?
\item Is it necessary to always pick the unexplored node with smallest optimal value to branch next, as we did above?
\item Similarly, do we always have to branch with the noninteger variable that is the furthest away from any integer value? 
\item If we were to branch on a node with value greater than the lower bound, 
would it be possible to immediately produce a better lower bound value?
\item If we were to branch on the only unexplored node with value matching the lower bound, 
would the next lower bound always come from a child of this node?
\end{enumerate}
\end{exercise}

\begin{exercise}
Build a branch-and-bound tree to solve the model from Exercise~\ref{ex:knapsack}. 
\end{exercise}

\begin{exercise}\label{ex:enumeration}
Back to the discussion about branch-and-bound and enumeration near Footnote~\ref{ft:enumeration}, 
there are two ways in which we avoid doing a complete enumeration with branch-and-bound.\footnote{If we were to use branch-and-bound to enumerate all solutions, 
we could reduce the computational cost of doing that by identifying nodes that would produce identical subtrees if explored~\cite{ebwlr}.} 
One of them is by fathoming unexplored nodes if their value implies that they would not lead to an optimal solution based on the upper bound. 
What is the other one? In other words, from where else could we have found other feasible solutions?
\end{exercise}

\begin{exercise}
There are many strategies that can be used to choose which node to explore next. 
On one extreme, we can try always branching on the deepest node available, 
which is known as Depth-First Search~(DFS). 
That is another way of interpreting what we did in the example along this section. 
On the other extreme, we can try always branching on the node that is closest to the root node, 
which is known as Breadth-First Search~(BFS).
If we had done that, 
we would have explored the two children of node 1 before exploring the children of those nodes. 
If we had explored the nodes 1, 2, and then 5 as in BFS, 
what would be the unexplored nodes after that as well as the upper and lower bounds of the model?
\end{exercise}

\section{Using a Mathematical Optimization Solver}\label{sec:listings}

Now we switch gears to discuss how to implement models on a computer. 
There are many code libraries---as well as software not requiring coding---that provide ways of expressing optimization models. 
They may either have their own optimization algorithms or interface with optimization solvers to get the model solved. 

We choose to use a computer programming language for greater flexibility; 
the Python language\footnote{If you learned more than one programming language, you may have been told that code ran in Python tends to be slower because it is often interpreted directly from the source code instead of compiled to a machine language executable. While that is true, the core implementation of the optimization algorithms in solvers is typically written in a lower level language such as C or C++ and then compiled, so what the Python libraries typically do is just interfacing with the solver executable. Moreover, the current popularity of Python and the abundance of libraries for manipulating data in Python make it the ideal language for a first interaction with an optimization solver.} due to its popularity; 
the Gurobi solver because it is regarded as one of the fastest solvers and, while being a commercial product, students with a \texttt{.edu} email address can get free academic licenses through their website~\cite{gurobi}; 
and the \texttt{gurobipy} library to have easier access to specific functions of this solver. 
If you have not used Python before, 
we recommend some free online books~\cite{py4e,boring}. 
For more about Gurobi, 
we recommend the online reference manual~\cite{gurobipy}.

We begin with a listing showing how to implement the model from Section~\ref{sec:intro}, 
which we denote as the Apple and Orange Model 1. We discuss it afterwards. \\

\begin{Verbatim}[numbers=left,xleftmargin=5mm,frame=single,label=Apple and Orange Model 1,fontsize=\small]
import gurobipy as gb

ao_model = gb.Model()

a = ao_model.addVar(vtype = gb.GRB.INTEGER)
o = ao_model.addVar(vtype = gb.GRB.INTEGER)

ao_model.setObjective(0.49*a + 0.50*o, gb.GRB.MINIMIZE)

ao_model.addConstr(214*a + 232*o >= 235)

ao_model.optimize()

print(ao_model.objVal)
print(a.X)
print(o.X)
\end{Verbatim}
The first line imports the \texttt{gurobipy} library and then renames it as \texttt{gb} to simplify the code.\footnote{You do not need to rename \texttt{gurobipy} or use \texttt{gb} if you do. However, \texttt{gb} is  a common choice and you may find examples using it online, similarly to the \texttt{pandas} library being called \texttt{pd}.} 
We create an object representing the model in Line 3, \texttt{ao\_model}.\footnote{Some people always use \texttt{model} for this variable. However, if you have multiple models, your code may end up referring to the wrong model---and detecting this kind of bug can be hard!} 
Then in Lines 5 and 6 we create two variables for that model. 
If we had not said that the variables were integer with \texttt{vtype = gb.GRB.INTEGER}, 
they would have been assumed to be continuous. 
By default, Gurobi assumes for each variable a lower bound of zero and no upper bound. 
That would be the same as proving the optional arguments \texttt{lb=0} and \texttt{ub = gb.GRB.INFINITY}. If we wanted no lower bound, we could similarly use \texttt{lb = - gb.GRB.INFINITY}.\footnote{Infinity has a value in Gurobi, although a big one, which is typically $10^{100}$ written as \texttt{1e+100}.}

Now it is a good moment to pause and talk about the variables \texttt{a} and \texttt{o} in the model. These are program variables representing the decision variables, but not integer variables as you would perhaps expect. They will have an integer value in them only after we have created and solved the model, but that is just one part of the information that they represent, which includes bounds and other things. When we multiply these variables or add them, we are creating an expression. Likewise, that expression can only be evaluated if we have values for the decision variables in it.

In Line 8 we provide an objective function to the model using the variables that we created. The second argument tells if the model is for minimization (with \texttt{gb.GRB.MINIMIZE}) or maximization (with \texttt{gb.GRB.MAXIMIZE}), and minimization is assumed if that argument is missing. In Line 10 we define the only constraint of the model, 
in which we write a linear expression and compare it with another value using \texttt{>=} for $\geq$. Other possible options are \texttt{==} for $=$ and \texttt{<=} for $\leq$.\footnote{Strict comparators such as $<$, $>$, and $\neq$ are typically not used in optimization models because the optimal solutions may not be defined. For example, if increasing some variable $x$ is beneficial to the objective function but we have $x < 2$ as a constraint, 
then we can keep finding slightly better values for $x$ that are not exactly 2. 
For example, 
we may start with $x = 1.9$, 
then move to $x = 1.99$, 
then $x = 1.999$, 
and so on. 
In many cases with integer variables we can replace strict comparisons with nonstrict ones. 
For example, $x < 2$ would be equivalent to $x \leq 1$ if $x$ can only be integer.}

Finally, the solver is called in Line 12. This is when all the output up to Line 25 of the next listing, Apple and Orangle Model Output, is produced. 
We will not see solver logs again, 
but they shown important information that you may need to know:
\begin{itemize}
\item Line 6 gives the problem dimensions. Nonzeros refer to any coefficient that is not zero in the objective function or in a constraint.
\item Line 8 gives the total number of decision variables of each type.
\item Line 14 tells if the solver was able to find a feasible solution before doing a rigorous search, such as through branch-and-bound. At this point, the solver uses a different family of methods known as heuristics, 
which we briefly discuss in Section~\ref{sec:heuristic}.
\item Lines 15 to 17 summarize the solver effort to simplify the model before solving it, 
which is called \emph{presolve}, and which may go as far as finding an optimal solution already---as it happened in this case. This is a very rich topic, for which \cite{presolve} gives an idea of the tricks used by optimization solvers.
\item Line 19 tells us how many branch-and-bound nodes were explored and for how long. 
Logs are longer and show frequent updates if branch-and-boud does happen.
\item Line 22 reports the number of solutions found and the value of each. 
\item Line 24 confirms that an optimal solution was found. If the problem is infeasible, meaning that no solution was found, that information would also be part of the log. If you assume that there is always going to be a solution and ignore the log, you may get into trouble! \\
\end{itemize} 

\begin{Verbatim}[numbers=left,xleftmargin=5mm,frame=single,label=Apple and Orange Model 1 Output,fontsize=\scriptsize]
Gurobi Optimizer version 11.0.3 build v11.0.3rc0 (win64 - Windows 11.0 (226

CPU model: Intel(R) Core(TM) Ultra 7 165H, instruction set [SSE2|AVX|AVX2]
Thread count: 16 physical cores, 22 logical processors, using up to 22 thre

Optimize a model with 1 rows, 2 columns and 2 nonzeros
Model fingerprint: 0x84338cfd
Variable types: 0 continuous, 2 integer (0 binary)
Coefficient statistics:
  Matrix range     [2e+02, 2e+02]
  Objective range  [5e-01, 5e-01]
  Bounds range     [0e+00, 0e+00]
  RHS range        [2e+02, 2e+02]
Found heuristic solution: objective 1.0000000
Presolve removed 1 rows and 2 columns
Presolve time: 0.00s
Presolve: All rows and columns removed

Explored 0 nodes (0 simplex iterations) in 0.01 seconds (0.00 work units)
Thread count was 1 (of 22 available processors)

Solution count 2: 0.98 1 

Optimal solution found (tolerance 1.00e-04)
Best objective 9.800000000000e-01, best bound 9.800000000000e-01, gap 0.000
0.98
2.0
-0.0
\end{Verbatim}

In the last lines of our first model, 
we print the optimal value and the value of each decision variable in the optimal solution found. 
Note that we access that value as an attribute, say \texttt{a.X} for variable \texttt{a}. 
The last decision variable has the puzzling value of minus zero, which is a common occurrence in numerical software, and it happens because one bit of the floating-point representation is reserved for the sign. In the particular case of the number zero, that bit is irrelevant. 

Speaking of numerical issues, Line 24 of the log mentions that the solution is optimal for a given tolerance. The smaller that tolerance, the slower it may take to solve a model. There are other parameters for numerical precision in solvers, such as the tolerance for a number closer to an integer value to be considered as integer. 
In other words, 
and thanks to the limited precision of floating-point arithmetic, 
welcome to the brave new world of approximately correct, optimal, and integer solutions!\footnote{There is also growing research in exact optimization~\cite{cook2013rational,eifler2023rational}, which as such tend to be slower.}

The next two listings show models implemented in a different way, 
especially for making the code shorter. 
These models follow better computer programming practices by separating the data, such as the price of fruit or their nutrition values, from the logic being implemented. 
In the second listing, Apple and Orange Model 2, 
we show another possibility for the the same mathematical model as before. 
The main change comes from Lines 9 and 10, 
in which we define \texttt{fruit} to be a collection of decision variables 
by using \texttt{.addVars} followed by the number of decision variables 
instead of \texttt{.addVar} for every decision variable separately. 
The variables will be represented as if in a Python list, 
with the first being \texttt{fruit[0]} and the second being \texttt{fruit[1]}, 
as we can see in the constraint definition in Lines 14 to 17 and when printing the solution at the end. 
In the call to create the decision variables we also provide the optional argument \texttt{obj}, 
which is a list with the coefficients for each decision variable in the objective function. 
Now we do not need to provide the objective function explicitly, 
but we play safe by specifying in Line 12 that this is a minimization model instead of assuming it by default. 
\\ 

\begin{Verbatim}[numbers=left,xleftmargin=5mm,frame=single,label=Apple and Orange Model 2,fontsize=\small]
import gurobipy as gb

fruit_price = [0.49, 0.50]
fruit_potassium = [214, 232]
potassium_need = 235

ao_model2 = gb.Model()

fruit = ao_model2.addVars(2, obj = fruit_price, 
                          vtype = gb.GRB.INTEGER)

ao_model2.ModelSense = gb.GRB.MINIMIZE

ao_model2.addConstr(
    fruit_potassium[0]*fruit[0] 
    + fruit_potassium[1]*fruit[1] 
        >= potassium_need)

ao_model2.optimize()

print(ao_model2.objVal)
print(fruit[0].X)
print(fruit[1].X)
\end{Verbatim}

In the third listing, Fruit Model (Incomplete), we show how to implement part of the model from Section~\ref{sec:model}.
Besides showing how to set an upper bound for the variables in Line 11, which is actually not needed in this case because we will not be using that much fruit, 
this listing contrasts two forms of specifying a constraint:
In Lines 16 to 21 we define the potassium constraint by explicitly listing each coefficient and decision variable. 
In Lines 23 to 27 we define the calcium constraint by using a summation on index $i$. 
The function \texttt{gb.quicksum} produces a summation that is roughly the same as writing $\sum\limits_{i=1}^4 \texttt{fruit\_calcium[i]}  * \texttt{fruit[i]} \geq \texttt{nutrient\_need[1]}$. 
This function works in a similar way as list compression in Python, 
with which we can create a list of integer numbers from 1 to 10 by using \texttt{[ i for i in range(1, 11)]}, for example. 
We leave the fiber constraint to the reader as Exercise~\ref{ex:fiber}, which is why the listing is incomplete.
The last differentiation of this listing is that we have the convenience of a \texttt{for} loop printing the value of the decision variables at the end. 
\\

\begin{Verbatim}[numbers=left,xleftmargin=5mm,frame=single,label=Fruit Model,fontsize=\small]
import gurobipy as gb

fruit_price = [0.49, 0.50, 0.51, 0.52]
fruit_potassium = [214, 232, 375, 155]
fruit_calcium = [12, 60.2, 5.75, 14.2]
fruit_fiber = [4.8, 2.8, 5.31, 5.52]
nutrient_need = [235, 65, 5.6]

fruit_model = gb.Model()

fruit = fruit_model.addVars(4, obj = fruit_price, ub = 10,
                            vtype = gb.GRB.INTEGER)

fruit_model.ModelSense = gb.GRB.MINIMIZE

fruit_model.addConstr(
    fruit_potassium[0]*fruit[0] 
    + fruit_potassium[1]*fruit[1] 
    + fruit_potassium[2]*fruit[2]
    + fruit_potassium[3]*fruit[3]
        >= nutrient_need[0])

fruit_model.addConstr(
    gb.quicksum(
        fruit_calcium[i]*fruit[i]
        for i in range(4))
        >= nutrient_need[1])

fruit_model.optimize()

print(fruit_model.objVal)
for i in range(4):
    print(fruit[i].X)
\end{Verbatim}

The listings above are meant to provide you with the basics to write simple models. 
There are many other elements that can be used when implementing models 
with \texttt{gurobipy}, 
and we will see a few more of those in future sections.

\begin{figure}[h]
\centering
\includegraphics[width=\textwidth]{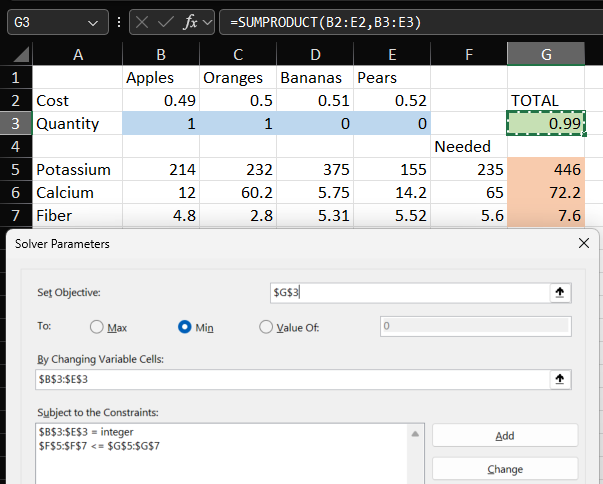}
\caption{Spreadsheet implementation of the ILP model from Section~\ref{sec:model} using Excel Solver.}
\label{fig:Excel}       
\end{figure}

\begin{svgraybox}
\textbf{Can I start modeling without coding?} \\
Many students start modeling optimization problems with spreadsheets, 
which is a common approach used by textbooks. 
Microsoft Excel has the built-in Solver, which needs to be activated as an Add-In, and that can solve problems with up to 200 binary variables. 
Figure~\ref{fig:Excel} shows the model from Section~\ref{sec:model} implemented with Excel Solver. 
Another free option that can be installed in Excel and also used as a plug-in in Google Sheets is OpenSolver~\cite{OpenSolver}, 
which uses the open source CBC solver from COIN-ON~\cite{coin-or}.
\end{svgraybox}

\begin{exercise}
Reimplement the potassium constraint with \texttt{gb.quicksum}.
\end{exercise}

\begin{exercise}\label{ex:fiber}
Implement the fiber constraint, which is missing from the last listing.
\end{exercise}

\begin{exercise}\label{ex:knapsackImpl}
Implement your model formulated for the knapsack problem from Exercise~\ref{ex:knapsack}.
\end{exercise}

\chproblem{\label{cp:randomKnapsack}
One of the goals of research is to help us better understand how something works. 
For example, we may want to know how hard it is to solve knapsack problems of a given size, or how much longer it takes to solve these problems as their size grow. 
By solving that problem once, for example with the size and coefficients from Exercise~\ref{ex:knapsack} that you may have implemented in Exercise~\ref{ex:knapsackImpl}, 
we are just getting started. 
Compared to other problems with similar size and coefficients as that particular knapsack problem, 
it is possible that the one that we implemented may have been solved significantly faster or slower than the average. 
In that sense, 
the runtime for that particular problem would be called \emph{anectodal evidence}, 
since it is a valid result but we do not know how it compares to the usual case. 

Follow the steps below to adapt your implementation from Exercise~\ref{ex:knapsackImpl} 
to solve a slightly different knapsack problem every time:
\begin{enumerate}[(a)]
\item Keep the number of items as 5 and the prices as they are in Table~\ref{tab:1}, 
but generate random integers between 1 and 10 for the happiness associated with each item. 
Run the code a few times to see how the solution and the time to solve change.\footnote{You can use the \texttt{random} module in Python for generating random integers. 
First, you need the line \texttt{import random} at the top of the code. Now you can use \texttt{random.randint(1,10)} at any point after importing the module to generate an integer number between 1 and 10. To generate a list of 5 values, you can use \texttt{[ random.randint(1,10) for i in range(5) ]} .}
\item Change the code to also generate random prices up to \$4.50 for each item.\footnote{You may potentially use \texttt{random.uniform(0.01,4.50)} for generating a floating point number between 0.01 and 4.50, but those numbers may have more than two decimal places. If you want to generate prices with an allowable fractional part, you can round that float output by using \texttt{round(random.uniform(0.01, 4.50), 2)} or generate a random integer for the total number of cents with \texttt{0.01*random.randint(1,450)}.}
\item Change the code to also make the number of items random between 5 and 20.
\end{enumerate}
}

\section{Assignment Problems and Graphs}\label{sec:assignment}

The discrete optimization models that we have discussed so far can be very useful when a quantity has to be integer, 
such as when we can only buy whole apples and oranges. 
However, 
some of the most interesting --- and challenging --- discrete optimization models have a different type of decision variable. 
In those models, 
we use binary variables for yes--no decisions. 
We have seen somes examples already: 
\begin{itemize}
\item In Exercise~\ref{ex:knapsack}, we had a knapsack problem in which we could only buy at most one bar of chocolate of each type. That was a smooth transition from integer to binary, since the value of the variable still represented a quantity.  
\item In Challenge Problem~\ref{ex:facility}, we had a facility location problem in which we had to decide which facilities to rent and, based on those, which facility would supply each store. 
The decision variables were provided to make it easier to understand how to model it. 
Still, this was a challenge problem because reasoning about constraints involving binary variables is different than thinking about quantities.\footnote{
Going back to the last item in part (a) of Challenge Problem~\ref{ex:facility}, 
you can use the constraint $x_{i j} \leq y_i$ to ensure that $x_{i j} = 1$ only if $y_i = 1$: 
If $y_i = 0$, then $x_{i j} = 0$. If $y_i = 1$, then $x_{i j}$ is not constrained.}
\end{itemize}

Now we will switch to problems that are very discrete in nature, 
in the sense that a solution for the LP relaxation is not as meaningful as 
it was for the previous models. 
We will start with some variants of the \emph{assignment problem} as a warm up to the next section, 
in which we will discuss the \emph{traveling salesperson problem}. 
That will help understanding common features of discrete optimization models. 

Let us suppose that a company has offices in the most populous cities in the United States, 
and that this company has some employees going from one office to another to audit the work of their colleagues. 
If every office has one of such employees and they all need to move to a different office every month, 
then we can try to minimize the distance that they need to travel. 

We will gradually develop this model across the rest of the section. 

\subsection{The First Model: One-Way Assignment}

To start modeling the problem, 
let us consider the fact that every employee needs to travel to another office. 
We start with three offices, numbered from 1 to 3, and we use $d_{i j}$ as the distance from $i$ to $j$.\footnote{
If $d_{i j} \neq d_{j i}$, 
we say that the problem is \emph{assymetric}. 
However, our example will be symmetric.
}  
That can be modeled as follows:
\[
\begin{aligned}
\min ~~& d_{1 2} x_{1 2} + d_{1 3} x_{1 3} + d_{2 1} x_{2 1} + d_{2 3} x_{2 3} + d_{3 1} x_{3 1} + d_{3 2} x_{3 2} \\
\text{s.t.} ~~&  x_{1 2} +x_{1 3} = 1 \\
& x_{2 1} +x_{2 3} = 1 \\
& x_{3 1} +x_{3 2} = 1 \\
& x_{1 2}, x_{1 3}, x_{2 1}, x_{2 3}, x_{3 1}, x_{3 2} \in \{0, 1\}
\end{aligned}
\]

We use a decision variable $x_{i j} \in \{0, 1\}$ to determine whether we go from $i$ to $j$, 
provided that $i \neq j$. 
In essence, this is an integer variable with lower bound of zero and upper bound of one.  
This is the most common type of discrete decision variable in MILP models\footnote{Quick reminder: the MILP models are a superset of the ILP models; hence more general.}, 
and it is sometimes called a binary or \emph{0--1} variable. 
It is very common to think about binary variables in logic terms as we do here, 
in which case the variable is either true (when it is 1) or false (when it is 0). 
In other words, $x_{i j} = 1$ means that we go from office $i$ to office $j$, with an impact of $d_{i j}$ on the total cost; if $x_{i j} = 0$, then we do not go directly from $i$ to $j$.

We present our implementation in listing One-Way Assignment Model:    
\\

\begin{Verbatim}[numbers=left,xleftmargin=5mm,frame=single,label=One-Way Assignment Model,fontsize=\small]
import gurobipy as gb

cities = ['New York-NY', 'Los Angeles-CA', 'Chicago-IL']

coordinates = {'New York-NY': (-74.0059413, 40.7127837),
 'Los Angeles-CA': (-118.2436849, 34.0522342),
 'Chicago-IL': (-87.6297982, 41.8781136)}

pairs = [('New York-NY', 'Los Angeles-CA'),
 ('New York-NY', 'Chicago-IL'),
 ('Los Angeles-CA', 'New York-NY'),
 ('Los Angeles-CA', 'Chicago-IL'),
 ('Chicago-IL', 'New York-NY'),
 ('Chicago-IL', 'Los Angeles-CA')]

owap_model = gb.Model()

x = owap_model.addVars(pairs, vtype = gb.GRB.BINARY)

owap_model.setObjective(
    2789.8*x[('New York-NY', 'Los Angeles-CA')]
    + 852.7*x[('New York-NY', 'Chicago-IL')]
    + 2789.8*x[('Los Angeles-CA', 'New York-NY')]
    + 1970.5*x[('Los Angeles-CA', 'Chicago-IL')]
    + 852.7*x[('Chicago-IL', 'New York-NY')]
    + 1970.5*x[('Chicago-IL', 'Los Angeles-CA')],
    gb.GRB.MINIMIZE
)

owap_model.addConstr( x[('New York-NY', 'Los Angeles-CA')]
                  + x[('New York-NY', 'Chicago-IL')] == 1)
owap_model.addConstr( x[('Los Angeles-CA', 'New York-NY')]
                  + x[('Los Angeles-CA', 'Chicago-IL')] == 1)
owap_model.addConstr( x[('Chicago-IL', 'New York-NY')]
                  + x[('Chicago-IL', 'Los Angeles-CA')] == 1)

owap_model.optimize()
visualize_solution(x, coordinates, pairs, distance)
\end{Verbatim}

One of the most important changes in terms of implementation 
is that we use \texttt{.addVars} to create a collection of decision variables that is not numbered sequentially from zero. 
Instead of that, we provide a list of keys for the decision variables, 
each of which consisting of a different pair of cities in which the offices are located,  
as if the variable \texttt{x} represented a dictionary in Python. 
For example, \texttt{x[('New York-NY', 'Los Angeles-CA')]} is the decision variable for 
whether the employee in the New York office goes to the Los Angeles office. 
We use these variables explicitly when defining the objective function and the first constraint, 
so that you can get used to them in the implementation of the model. 

Another impotant change is the call to the function \texttt{visualize\_solution} in the last line. 
There is a listing in the appendix with the code of this and other functions, 
which later will also be used to automatically produce the names of the cities and the distances between them by reading their names and coordinates from a file 
instead of having them as part of the code as in the listing above. 
From calling that function, 
we obtain an illustration of the optimal solution found in Figure~\ref{fig:TSP01}, 
which is discussed in Section~\ref{sec:two}. 

Now is good point to talk about what makes models with binary decision variables different: 
the LP relaxation tends to be less meaningful. 
For example, 
what do we make about going 0.75 from New York to Chicago and 0.25 from New York to Los Angeles? 
This is a little more convoluted than talking about slicing fruit.

\subsection{Modeling with Graphs}

\begin{figure}[h]
\centering
\includegraphics[width=\textwidth]{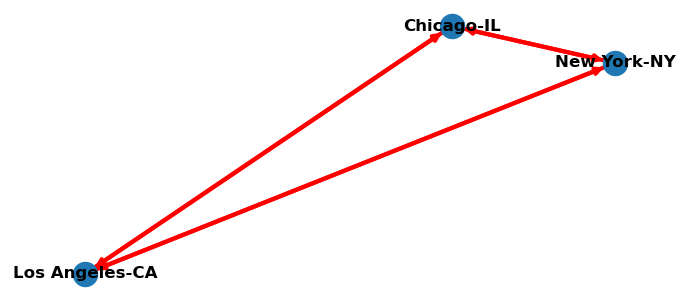}
\caption{Complete directed graph with the three most populous cities in the United States.}
\label{fig:TSP00}       
\end{figure}

Note that it does not make sense to go from one office to the same office in our case, 
and it could also be that there is no way to go directly from one office to another. 
To make our model more general but still compact to represent, 
we will use a \emph{graph}. 
By graph we mean a mathematical object that can be tought of as a set of dots connected by arrows (if \emph{directed}, representing a single direction) or line segments (if \emph{undirected}, representing both directions). 
Because graphs are very useful modeling devices, 
there is a subject dedicated to the study of graphs that is known as \emph{graph theory}~\cite{murty}.  
In our case, 
the dots are the offices and the arrows represent if we can travel between them. 
We will use arrows because the decision about an employee going from office $i$ to office $j$ ($x_{i j}$) is not the same as the decision about an employee going from office $j$ to office $i$ ($x_{j i}$).  
In Figure~\ref{fig:TSP00} we show a directed graph with the three most populous cities in the United States: New York, Los Angeles, and Chicago. 
For each pair of cities, there are two opposite arrows overlapping. 

More formally, 
we can represent a directed graph as $G = (N, A)$, 
where $N := \{ 1, \ldots, n\}$ is the set of nodes --- which represent offices numbered from 1 to $n$ in our case ---  
and $A \subseteq \{ (i, j) : i \in N, j \in N \setminus \{ i \} \}$ is the set of arcs connecting different nodes, 
possibly with some omissions --- if we cannot go from one city to another.  
We say that a graph is \emph{complete} is all nodes are connected to all other nodes.

To rewrite the optimization model using a graph, 
we need to refer to what offices can be reached directly from a given office and to what offices can reach directly a given office. 
Following standard convention in graph theory, 
we will use $\delta^+(i) \subseteq \{ j : j \in N \setminus \{ i \} \}$ to refer to all the office that can be visited immediately after office $i$, 
and $\delta^-(j) \subseteq \{ i : i \in N \setminus \{ j \} \}$ to refer to all the office that can be visited immediately before city $j$.\footnote{It is very common to have a problem with two indices, say $i$ and $j$, and it tends to be helpful not to flip them in any expression even if they refer to the same type of entity--- this is why we talk about $\delta^+(i)$ if departing office $i$ but then we talk about $\delta^-(j)$ if arriving at office $j$.} 
Now we can generalize the model as follows:

\[
\begin{aligned}
\min ~~& \sum_{(i,j) \in A} d_{i j} x_{i j} \\
\text{s.t.} ~~&  \sum_{(i,j) \in \delta^+(i)} x_{i j} = 1 & \forall i \in N \\
& \textbf{x} \in \{0, 1\}^{|A|}
\end{aligned}
\]

We use the symbol $\forall$, which is known as \emph{for all}, 
to express that there is a constraint associated with each office $i \in N$. 
We also use summations over sets, such as $\sum\limits_{(i,j) \in A}$ to express that we are summing over all the arcs in the graph.  

\subsection{The Second Model: Two-Way Assignment}\label{sec:two}

\begin{figure}[h]
\centering
\includegraphics[width=\textwidth]{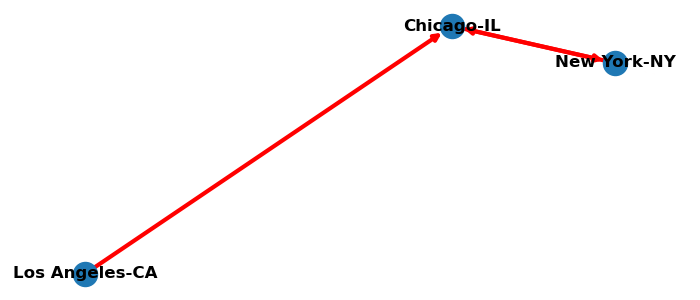}
\caption{Optimal solution of the One-Way Assignment Model on the three most populous cities in the United States: from Los Angeles to Chicago and from Chicago and New York to each other.}
\label{fig:TSP01}       
\end{figure}

Figure~\ref{fig:TSP01} illustrates the optimal solution obtained from listing One-Way Assignment Model.  
We can see that 
every employee is going from one office to another. 
However, not every office is receiving a new employee. 
Namely, 
the employee in Los Angeles goes to Chicago, 
the employee in Chicago goes to New York, 
but the employe in New York also goes to Chicago --- and nobody goes to Los Angeles. 
Hence, we are left with two employees in Chicago and none in Los Angeles. 

We add a second set of constraints to the model as shown below, which looks very similar to the first set of constraints that already existed in the previous model: 
\[
\begin{aligned}
\min ~~& \sum_{(i,j) \in A} d_{i j} x_{i j} \\
\text{s.t.} ~~&  \sum_{(i,j) \in \delta^+(i)} x_{i j} = 1 & \forall i \in N \\
& \sum_{(i,j) \in \delta^-(j)} x_{i j} = 1 & \forall j \in N \\
& \textbf{x} \in \{0, 1\}^{|A|}
\end{aligned}
\]

With the first constraint that was already present in the previous model, we fix the office of origin $i$ and sum over all the destinations with $j$ to ensure that the employee in office $i$ goes somewhere else. 
With the new constraint that was introduced in this model, we fix the office of destination $j$ and sum over all the destinations with $i$ to ensure that every office $j$ receives an employee from somewhere else. 

We add this new type of constraint and make the implementation more abstract 
by not including coefficients directly in the next listing, Two-Way Assignment Model:
\\

\begin{Verbatim}[numbers=left,xleftmargin=5mm,frame=single,label=Two-Way Assignment Model,fontsize=\small]
import gurobipy as gb

cities, coordinates, distance = read_problem(
    "1000-cities.csv", 5)

twap_model = gb.Model()

x = twap_model.addVars(distance, obj = distance, 
		           vtype = gb.GRB.BINARY)

for i in cities:
    twap_model.addConstr(
        gb.quicksum(
            x[i,j] for j in cities 
                   if j!=i)
        == 1)

for j in cities:
    twap_model.addConstr(
        gb.quicksum(
            x[i,j] for i in cities 
                   if j!=i)
        == 1)

twap_model.optimize()
visualize_solution(x, coordinates, distance)
\end{Verbatim}

We can now choose to work with up to 1000 cities instead of only 3.
All the data that was explicitly represented in the previous listing 
is now loaded and processed based on the text file \texttt{1000-cities.csv}. 
The first 20 lines of that file are available in a listing in the appendix, 
whereas the entire file is available at \url{https://github.com/thserra/discreet} .
We use the auxiliary function \texttt{read\_problem} 
from a listing in the appendix to read and process that text file.

In Line 4, 
we load the data to solve the problem with $n = 5$ cities. 
We can change this line to another integer number up to $n = 1000$ to solve a problem of different size. 
In Line 8, we assign the coefficients of the objective function using the new dictionary \texttt{distance}, 
which maps each pair of cities to their corresponding distance. 
Since we are not changing \texttt{twap.ModelSense} as we did with \texttt{fruit\_model.ModelSense} in the last listings of Section~\ref{sec:listings}, 
we are implicitly assuming that Gurobi will minimize the objective function. 
In Lines 11 to 16 we implement the first constraint with a \texttt{for} loop on $i$, 
where the constraint in each city $i$ is expressed with a summation of every city $j$ other than $i$. 
Note that this is less general than the ILP model: 
for simplicity of presentation, 
we assume in this implementation that graph $G$ is complete, 
meaning that from any city it is possible to directly visit any other city. 
In Lines 18 to 23 we implement the second constraint in a similar way, 
but looping instead on $j$ and varying $i$ within the summation. 

\begin{figure}[h]
\centering
\includegraphics[width=\textwidth]{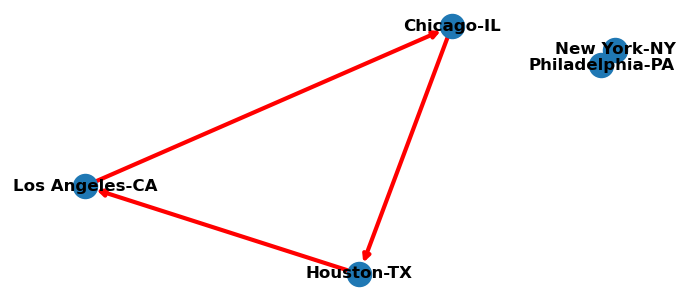}
\caption{Optimal solution of the Two-Way Assignment Model on the five most populous cities in the United States: one subtour has New York and Philadelphia; the other has Chicago, Houston, and Los Angeles.}
\label{fig:TSP02}       
\end{figure}

Figure~\ref{fig:TSP02} illustrates an optimal solution for the Two-Way Assignment Model with offices in five cities. 
In this solution, 
two employees alternate between New York and Philadelphia, 
whereas another three employees cycle among Chicago, Houston, and Los Angeles. 
Each of those sets of cities is defining what we call a \emph{subtour}, 
meaning that there is a sequence of arcs connecting those cities and such that we enter and leave each city once. 
If there was a single set with all the cities, 
that would be a \emph{tour}. 
Since the first constraint ensures that one arc leaves every city and the second constraint ensures that one arc enters every city,
then by starting in any given city we will visit a sequence of cities without any repetition (since a repetition would imply two arcs entering a same city) up to the point where we return to the first city (since we have not entered that city yet, and we may eventually have no other city left to visit). 
Hence, 
the solution will always have a tour or multiple subtours.

\begin{exercise}
Answer to the following questions:
\begin{enumerate}[(a)]
\item Is the optimal value of One-Way Assignment Model greater, smaller, or equal to the optimal value of Two-Way Assignment Model?
\item What about the optimal values of their LP relaxations?
\item Would you expect that to always be the case? If not, under which conditions would this keep happening?
\end{enumerate}
\end{exercise}

\begin{exercise}
One common type of assignment problem is to assign people to do tasks. 
Each person may take a different amount of time to finish each of the tasks. 
Each task must be assigned to a person. 
\begin{enumerate}[(a)]
\item Create an example in which you have 5 people and 5 tasks, 
and your goal is to find an assignment in which the total time to complete the tasks is minimized. 
\item Implement this model in Python and solve it.\footnote{You may find it convenient to have a matrix for how long each person takes to complete each task. To create a $2 \times 3$ matrix in Python, you can use $M = [ [ 5 ~ 6 ~ 7 ], ~ [ 2 ~ 9 ~ 3 ]]$. Then you can use $M[0][2]$ to access the element in the first row (numbered 0) and third column (numbered 2).}
\item Now remove the last task from your example and change the model to assign one person to each task. Note that one person will not be assigned to any task.
\end{enumerate}
\end{exercise}

\chproblem{\label{cp:randomFacility}
Implement your model formulated for the facility location problem in Challenge Problem~\ref{ex:facility} with random generation for the problem size and coefficients, similar to what you have done in Challenge Problem~\ref{cp:randomKnapsack}. 

In order to generate consistent distances between the locations, 
you may consider generating the $x$ and $y$ coordinates of each location as random floating point numbers from 0 to 10, for example. 
Then you may refer to the listings presented in this section for calculating the distance between locations.\footnote{You may also find it convenient to use a dictionary of tuples for the coordinates of each location. To represent the coordinates of two locations in Python numbered 0 and 1, you can use \texttt{locations = \{ 0 : (1.5, 3.2), 1 : (7.5, 3.2) \}}. Then you can use \texttt{locations[1][0]} to access the first coordinate (numbered 0) of location 1.}
}

\section{The Traveling Salesperson Problem}\label{sec:tsp}

Going back to the example from the previous section, 
what if the company expects that every employee rotates across all the offices? 
That changes it from a simpler assignment problem 
to the classic Traveling Salesperson Problem~(TSP). 
We are now crossing the line between the problems that are easy and hard to solve.\footnote{In general, finding an optimal solution for an ILP model is a hard problem. That means that trying to solve it by using branch-and-bound may sometimes produce an exponential number of explored nodes. But there are special cases, such as the assignment problem, in which we know specialized algorithms that are guaranteed to solve the problem efficiently. Despite that, using a general-purpose optimization solver for this and other easy problems often works well in practice.}

The Traveling Salesperson Problem (TSP) 
is commonly used as one of the first discrete optimization problems when introducing the topic because it is easy to understand. 
This is also a problem that quickly gets harder to solve as it grows in size if we are not careful. 
But even if we are careful, it eventually gets too hard anyway!

In this problem, 
we have $n$ cities that we would like to visit, numbered from 1 to $n$, 
and we would like to minimize the total distance traveled. 
The solution should be a \emph{tour}, which visits each city once and then returns to the city from where we started.  
Note that any city can be chosen as the first city, 
since we must visit all the cities. 
Moreover, 
we can consider the models from the last section as a starting point. 
The solution of each of those models made it evident that other constraints are needed, 
which leads to a third version of that same model. 

To proceed, 
we need to ``break'' those subtours in some way. 
For example, by having a constraint that prevents a subtour involving Chicago, Houston, and Los Angeles. 
Let us say that we number those cities as 1, 2, and 3 --- and thus have New York and Philadelphia as 4 and 5. 
In this case, we can use the constraint
\[
x_{1 4} + x_{1 5} + x_{2 4} + x_{2 5} + x_{3 4} + x_{3 5} \geq 1
\]
to ensure that we go at least once from a city in the first group (Chicago, Houston, and Los Angeles) 
to a city in the second group (New York and Philadelphia).  
This is known as a \emph{subtour elimination constraint}: 
If this constraint is not satisfied, 
then we would have no arc from a city in the first group to a city in the second group. 
Hence, there would be no way to reach New York or Philadelphia from Chicago, Houston, or Los Angeles; 
and there would be subtours involving the cities in each group. 

Next we will consider how to write the same constraint for any subset of nodes. 
Let us assume for now that the graph is complete. 
For a group of cities associated with a set of nodes $S$, 
which in the case above would be $S = \{1, 2, 3\}$, 
we can use 
\[
\sum_{i \in S} \sum_{j \in N \setminus S} x_{i j} \geq 1.
\]
to ensure that we always depart from at least one node in the subset $S$ to another node not in $S$ (hence in $N \setminus S$)\footnote{In set theory notation, $N \setminus S$ is the set obtained by subtracting from the set $N$ the elements that are also in set $S$.} in any feasible solution. 

Now let us consider that the graph is not complete. 
For every node $i$ considered in the outer  summation, 
the set of nodes that can be directly reached is given by $\delta^+(i)$. 
Hence, the set of nodes for inner summation depends on $i$. So we have
\[
\sum_{i \in S} \sum_{j \in N \setminus S} x_{i j} \geq 1
\]
as the general form for one subtour elimination constraint. 

Note that we could have expressed the same restriction using other constraints. 
We will discuss in more detail in the next section 
why that would not be a good as the option above. 
Such a discussion is important because it comes with a form of analysis that 
may help you think about how to better formulate new problems. 

\begin{figure}[h]
\centering
\includegraphics[width=\textwidth]{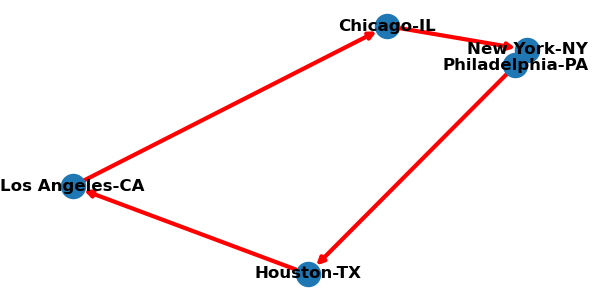}
\caption{Optimal solution for the TSP on the five most populous cities in the United States: from New York to Philadelphia, then Houston, Los Angeles, Chicago, and New York again.}
\label{fig:TSP03}       
\end{figure}

By using a subtour elimination constraint for every proper subset $S$ of the set of nodes $N$ with size at least 2, 
we obtain the following TSP model:  
\[
\begin{aligned}
\min ~~& \sum_{(i,j) \in A} d_{i j} x_{i j} \\
\text{s.t.} ~~&  \sum_{(i,j) \in \delta^+(i)} x_{i j} = 1 & \forall i \in N \\
& \sum_{(i,j) \in \delta^-(j)} x_{i j} = 1 & \forall j \in N \\
& \sum_{i \in S} \sum_{j \in \delta^+(i) \setminus S} x_{i j} \geq 1 & \forall S \subset N, |S| \geq 2 \\
& \textbf{x} \in \{0, 1\}^{|A|}
\end{aligned}
\]
For the same \emph{instance}\footnote{In optimization we often use the word problem to talk about an optimization problem in the abstract, 
and we use the word instance for a specific case. 
When we write a TSP model in which the number of cities is given by $n$ 
and the distance between a pair of cities is $d_{i j}$, 
this model can be used to solve instances in which the values of $n$ and $d_{i j}$ vary. 
In that sense, 
the listing One-Way Assignment Model is written with one instance in mind (all values in the code) 
and the listing Two-Way Assignment Model is written with the problem in mind (most of the data is loaded from a file). } of five cities for which we have seen the optimal solution for the Two-Way Assignment Model in shown Figure~\ref{fig:TSP02},
we now have an optimal solution satisfying all the subtour elimination constraints in Figure~\ref{fig:TSP03}.

We present next the listing TSP Model, with all the constraints needed:
\\

\begin{Verbatim}[numbers=left,xleftmargin=5mm,frame=single,label=TSP Model,fontsize=\small]
import gurobipy as gb
import time 

cities, coordinates, distance = read_problem(
    "1000-cities.csv", 5)

tsp_model = gb.Model()

x = tsp_model.addVars(distance, obj = distance, 
		           vtype = gb.GRB.BINARY)

for i in cities:
    tsp_model.addConstr(
        gb.quicksum(
            x[i,j] for j in cities 
                   if j!=i)
        == 1)

for j in cities:
    tsp_model.addConstr(
        gb.quicksum(
            x[i,j] for i in cities 
                   if j!=i)
        == 1)

tsp_model.addConstrs(
    gb.quicksum( 
        x[i,j]
        for i in subtour
        for j in cities if j not in subtour) 
    >= 1
    for subtour in subtours(cities)
)

start_time = time.time()
tsp_model.optimize()
end_time = time.time()
visualize_solution(x, coordinates, distance)
print("Duration:", end_time - start_time)
\end{Verbatim}
When expressing the new type of constraint in Lines 26 to 33, 
we use the method \texttt{.addConstrs}. 
The difference between \texttt{.addConstr} and \texttt{.addConstrs} is, in a sense, 
similar to the different between \texttt{.addVar} and \texttt{.addVars}: 
we use a single instruction to define multiple constraints with \texttt{.addConstrs}. 
You may also note in Line 32 that we use the Python syntax for list comprehension 
for indexing the constraints. 
We have also been using it for summations within each \texttt{gb.quicksum}, 
as in Lines 29 and 30. 
Hence, 
it is as if we were looping over all the subsets of cities (with Line 32) 
and for each subset we sum up the decision variables for the arcs leaving the subset (with Lines 27 to 30) 
in order to create the constraint preventing the subtours involving that subset. 
Therefore, anything that we do with \texttt{.addConstrs} can also be done with external loops over \texttt{.addConstr} and vice versa. 

Another difference in this listing 
is that we measure the time before and after calling the solver in Line 36. 
This allows us to measure how many seconds it takes to solve the model, 
and then compare the runtimes when we implement the subtour elimination constraints in a different way in Section~\ref{sec:callback}.

\begin{exercise}
Let us compare lower bounds on the optimal value of the TSP: 
\begin{enumerate}[(a)]
\item For the instance with 5 cities used in TSP Model,  
adapt the code to compute the following:
\begin{enumerate}[(i)]
\item the optimal value of the ILP model with the first type of constraint only;
\item the optimal value of the LP relaxation with the first type of constraint only;
\item the optimal value of the ILP model with the first two types of constraints only;
\item the optimal value of the LP relaxation with the first two types of constraints only; and
\item the optimal value of the LP relaxation of the model with the three types of constraints. 
\end{enumerate}
\item What comparisons among those values would always hold? For example, is it true that (i) $\geq$ (ii) for any instance?
\end{enumerate}
\end{exercise}

\begin{exercise}
The right-hand side value of 1 in a subtour elimination constraint over a subset of nodes $S$
means that we have at least one arc from a node in $S$ to a node not in $S$. 
If we had used an equality instead, 
we would have exactly one arc from a node in $S$ to a node not in $S$, 
but that remove feasible solutions. 
For the tour represented in Figure~\ref{fig:TSP03}, 
what subset of nodes $S$ would have more than one departing arc?
\end{exercise}

\begin{exercise}
The subtour elimination constraints ensure that we leave the nodes in each subset at least once. 
Why don't we need a constraint ensuring that we also enter the nodes in each subset at least once?
\end{exercise}

\begin{exercise}
Do we need to break subtours of all sizes? If not, what is the largest size for which we need to break a subtour, in terms of the number of cities $n$?
\end{exercise}

\begin{exercise}\label{ex:addConstrs}
Create two new versions of the code in listing TSP Model, 
where (1) only the method \texttt{.addConstr} is used for defining constraints; 
and (2) only the method \texttt{.addConstrs} is used instead.
\end{exercise}

\chproblem{\label{cp:knapsackStudy}
Now that we seen how to measure the runtime of a model, 
we may conduct a study of the time that it takes to solve a problem according to its size. 
We will start with the knapsack problem model that you may have implemented in Challenge Problem~\ref{cp:randomKnapsack}. 

Follow the steps below to adapt your implementation from Challenge Problem~\ref{cp:randomKnapsack} to produce a computational experiment to evaluate the time that it takes to solves knapsack problems according to their size:
\begin{enumerate}[(a)]
\item Encapsulate all the code from Challenge Problem~\ref{cp:randomKnapsack} in a function that takes the number of items of the knapsack as an argument. All other coefficients should be randomly generated as before. Surround the \text{.optimize} call with instructions for measuring the duration of the call. Make the function return that duration.
\item Call the function that you created at least 30 times with the same number of items, say 20. 
Save all the values returned to a list. 
Plot a boxplot and a histogram of those runtimes using your favorite statistical visualization tool.
\item Now let us create a loop for testing different number of items, say from 10 to 100 in increments of 10. 
For each of those cases, call the function at least 30 times. 
Now plot side by side the boxplot of runtimes for each number of items. 
\item Draw your conclusions about how the number of items affects the median runtime and the variability of the runtimes. 
\item Can you formulate a hypothesis? If so, conduct more experiments with problems having even more items to test if your hypothesis holds. 
\end{enumerate}
}

\chproblem{
Follow similar steps as in Challenge Problem~\ref{cp:knapsackStudy} with the facility location problem model that you may have implemented for Challenge Problem~\ref{cp:randomFacility}. 
Conduct a study by solving problems using different numbers of facilities and stores. 
Choose at least five values for each of those.  
For each combination, run between 30 and 50 experiments. 
Analyze the results to try finding patterns by, 
for example, 
producing a matrix of average runtimes according to the numbers of facilities and stores; 
or producing a matrix of boxplots for each of those combinations. 
}

\section{Modeling Tricks and Model Strength}\label{sec:tricks}

When we come across a new problem, 
we may need some time to think about how to formulate a linear model for that problem---if at all possible. 
However, there are some tricks that may help us come up with the constraints that we need more easily, 
in particular when the decision variables are binary. 
Moreover, since we may realize that more than one formulation is possible for an optimization problem, 
there are also ways to compare modeling choices to determine which one is preferable to use. 

For one reason or another, 
we may sometimes need to remove one binary solution, 
or perhaps prevent that a subset of binary variables take a specific assignment. 
For example, if we were to remove $(\bar{x}_1, \bar{x}_2, \bar{x}_3) = (0, 0, 0)$, 
then we could use either    
\[
x_1 + x_2 + x_3 \geq 1
\]
or 
\[
x_1 + x_2 + x_3 \geq 0.3 .
\]
If we consider that variables $x_1$, $x_2$, and $x_3$ can only take binary values, 
both of those constraints would only remove $(x_1, x_2, x_3) = (0, 0, 0)$. 
However, 
if we consider that variables $x_1$, $x_2$, and $x_3$ could also take noninteger (also called \emph{fractional}) values, 
such as when solving the LP relaxation of an ILP model, 
then those constraints are not equivalent: 
Some solutions removed by the first constraint would not be removed by the second constraint, 
such as $(x_1, x_2, x_3) = (0.1, 0.1, 0.1)$, 
whereas any solution removed by the second constraint is also removed by the first since the first constraint ($x1 + x_2 + x_3 \geq 1$) implies the second constraint ($x_1 + x_2  + x_3 \geq 0.3$). 

If we were to add the first constraint to an ILP model that previously had $(x_1, x_2, x_3) = (0, 0, 0)$ as a feasible solution, 
then we may take less time to solve the resulting model than if we had added the second constraint instead: 
In the first case we would have fewer fractional solutions, 
and thus possibly fewer nodes of the branch-and-bound tree in which the optimal solution of the LP relaxation is fractional. 
Hence, we may solve the problem faster with the first constraint.  

There are some ways by which we can make a constraint with $\geq$ stronger, 
in the sense that it removes more fractional solutions. 
First, we can decrease a coefficient on the Left-Hand Side~(LHS), 
i.e., we decrease by how much we multiply one of decision variables. 
Second, we can increase the coefficient on the Right-Hand Side~(RHS), 
i.e., we increase the constant not multiplying any decision variable. 
However, 
we should make sure that we are not removing any other binary solutions. 

Because a linear constraint removes a continuous set of solutions, 
verifying if we are removing other integer solutions int his case only requires verifying if we are removing any neighbors of the solution that we want to remove. 
In the case of $(x_1, x_2, x_3) = (0, 0, 0)$, 
those would be $(1, 0, 0), (0, 1, 0)$, and $(0, 0, 1)$. 
In fact, 
we can make the second constraint stronger with either of the options above: 
$0.5 x_1 + x_2 + x_3 \geq 0.3$ and $x_1 + x_2 + x_3 \geq 0.5$ would both remove $(x_1, x_2, x_3) = (0.1, 0.1, 0.1)$ while not removing any of the neighbors. 
However, doing that with the first constraint would invariably remove at least one of the neighbors. 
For the first constraint, 
note that the LHS value for all the neighbors, 
in this case the expression $x_1 + x_2 + x_3$,
is equal to the RHS, 
in this case the constant 1. 
When those two values match for a given solution, 
we say that the constraint is \emph{active} on the solution. 
That prevents us from strengthening the constraint, 
which is actually a good sign.

More generally, let $\textbf{x} = \bar{\textbf{x}}$ be such a solution on $n$ binary variables that we would like to remove from a model. 
We can express a constraint to remove $\textbf{x} = \bar{\textbf{x}}$ as follows:
\[
\sum_{i \in \{1, \ldots, n\} : \bar{x}_i  = 1} x_i + \sum_{i \in \{1, \ldots, n\} : \bar{x}_i  = 0} (1 - x_i) \leq n - 1.
\]
This is known as the canonical cut on the unit hypercube~\cite{balas1972cuts}. 
If you think about the binary solutions on $n$ variables as representing the vertices of an $n$-dimensional hypercube (i.e., a line segment for $n=1$, a square for $n=2$, a cube for $n=3$, and so on), 
then this constraint removes exactly one vertex along with as much of the fractional solutions as possible. 
In other words, 
this constraint is active on all the neighbors of the solution being removed. 
For $(\bar{x}_1, \bar{x}_2, \bar{x}_3) = (0, 0, 0)$, we would have 
\[
(1 - x_1) + (1 - x_2) + (1 - x_3) \leq 2,
\]
which is equivalent to that first constraint ($x_1 + x_2 + x_3 \geq 1$). 

We can use such cuts for removing an optimal solution of the model and then solving the model again. 
Doing that may produce another optimal solution, or the second best solution instead. 
This can be useful when there is uncertainty on the data, 
in which case it is not entirely clear if the first solution is indeed the best option; 
or if the first solution obtained is not easy to implement, 
which may happen due to aspects of the real-world problem that were not accounted for in the model~\cite{camm2014asp}.

In practice, we say that one model is stronger than another for the same problem if it has fewer fractional solutions. 
As seen in Chapter~\ref{sec:bb}, 
branch-and-bound can be a lengthy process if too many LP relaxations have solutions with fractional values for the integer variables of the ILP model, 
to the point of taking an exponential time on the size of the instance.  
Note that using canonical cuts on the unit hypercut is a first step, 
but it may not be sufficient: 
It is possible to find a constraint that removes even more fractional solutions if one of the neighbors of $\bar{\textbf{x}}$ is infeasible. 
More generally, we need to compare entire ILP models for the same problem, 
in the hope of finding which ones have the smallest LP relaxations.\footnote{
Whereas solutions on integer variables are countable and finite if the decision variables are bounded, 
solutions on continuous variables are not. Let's think of those as a continuous set.
}
To focus on models that we can compare directly, 
consider the definitions below: 
\begin{definition}
Two models are \emph{equivalent} if they have the same feasible set.
\end{definition}

\begin{definition}
For two equivalent ILP models $A$ and $B$, 
with $R(A)$ and $R(B)$ as the feasible set of their LP relaxations, 
model $A$ \emph{dominates} $B$ if $R(A) \subseteq R(B)$. 
\end{definition}

\begin{definition}
For two equivalent ILP models $A$ and $B$, 
with $R(A)$ and $R(B)$ as the feasible set of their LP relaxations, 
model $A$ \emph{strictly dominates} $B$ if $R(A) \subset R(B)$. 
\end{definition}

In other words, 
all ILP models for the same problem on the same decision variables are equivalent. 
These models may dominate each other, in case their LP relaxations are equivalent; or 
one of them may strictly dominate the other; or 
no dominance relaxation may exist between them if each LP relaxation has a solution that the other does not have. 
When comparing models, our main interest is on strict dominance. 
We say that a model is \emph{stronger} if it strictly dominates the other model. 

For an example of how to compare two forms of modeling the same requirement, 
let us consider again the removal of subtours in the TSP. 
We will use the same setting as in the previous section, 
with $n=5$  nodes and the subset of nodes $S = \{1, 2, 3\}$. 
We can avoid the subtour $1 \rightarrow 2 \rightarrow 3 \rightarrow 1$ with the canonical cut 
\[
x_{1 2} + x_{2 3} + x_{3 1} \leq 2. 
\]
However, 
there is another subtour involving $S$, $1 \rightarrow 3 \rightarrow 2 \rightarrow 1$, 
leading to the cut 
\[
x_{1 3} + x_{3 2} + x_{2 1} \leq 2.
\]
How would these two cuts together compare with the subtour elimination constraint 
\[
x_{1 4} + x_{1 5} + x_{2 4} + x_{2 5} + x_{3 4} + x_{3 5} \geq 1
\]
that we have seen before?

First, note that this is not about the number of constraints involved: 
One model may strictly dominate another while having less or more constraints in total. 
Second, those constraints as presented are defined in different spaces. 
Hence, we need to put them in terms of the same decision variables for comparison. 

On the one hand, given the first set of constraints in TSP Model, 
\[
\sum_{j \in \{1, \ldots, 5\} \setminus \{ i \}} x_{i j} = 1 \qquad \forall i \in \{1, \ldots 5\}, 
\]
we know that $x_{1 4} + x_{1 5} = 1 - x_{1 2} - x_{1 3}$, 
$x_{2 4} + x_{2 5} = 1 - x_{2 1} - x_{2 3}$, and 
$x_{3 4} + x_{3 5} = 1 - x_{3 1} - x_{3 2}$. 
Hence, the subtour elimination constraint can be restated as 
\[
(1 - x_{1 2} - x_{1 3}) + (1 - x_{2 1} - x_{2 3}) + (1 - x_{3 1} - x_{3 2}) \geq 1, \text{ or}
\] 
\[
(x_{1 2} + x_{2 3} + x_{3 1}) + (x_{1 3} + x_{3 2} + x_{2 1}) \leq 2.
\]
Since $\textbf{x} \geq 0$, 
then the subtour elimination constraint implies both cuts used for eliminating specific subtours. 
On the other hand, 
neither of the cuts removing subtours would remove  
$(x_{1 2}, x_{1 3}, x_{2 1}, x_{2 3} , x_{3 1}, x_{3 2}) = (0.5, 0.5, 0.5, 0.5, 0.5, 0.5)$, 
which is a fractional solution removed by the classic subtour elimination constraint.
Therefore, 
the TSP Model with the classic subtour elimination constraints strictly dominates an alternate model with cuts removing each subtour. 
Better yet, the stronger model also has fewer constraints!\footnote{
Although model strength is important, 
an excessive number of constraints may certainly slow down the time for solving the LP relaxation at each node of the branch-and-bound tree. 
According to Bill Cook, sometimes solving the LP relaxation is the bottleneck in large TSP models. 
That may also be an issue in cases where we use constraints to reduce the feasible set when there is symmetry in the model. 
These are called symmetry-breaking constraints~\cite{Margot2010} and they often help solving the problem faster, 
but in some cases the number of such constraints is so large that adding them to the model increases the runtime. 
One alternative in those cases is to only exploit the symmetries by building a different type of branch-and-bound tree, 
such as with orbital branching~\cite{orbital}.  
}

\begin{exercise}
Using the TSP model with $n = 5$, generate and remove optimal solutions until you get a solution worse than the optimal value. How many solutions did you get? Are they entirely distinct from one another?
\end{exercise}

\begin{exercise}
In production scheduling problems, it may be important to capture when transitions happen in machines, 
since they may imply a setup  cost --- for example, if cleaning is required. 
Consider the following setting:
\begin{itemize}
\item $x$ is a binary variable denoting if today a given production line is producing bottles of soda;  
\item $y$ is a binary variable denoting if tomorrow that same production line is producing mineral water bottles; and 
\item  and $z$ is a binary variable that is 1 if, and only if, $x=1$ and $y=1$. 
\end{itemize}
How can we write constraints for modeling that?
\end{exercise}

\begin{exercise}
Suppose that $x, y,$ and $z$ are binary decision variables, 
and that we want to have them in a model along with a variable $w = x y z$ (i.e, $w$ is equal to the other three multiplied).  
Is it possible to formulate that as a linear model?
\end{exercise}

\begin{exercise}
We must decide which lots to buy for developing a real estate complex. There are four consecutive lots, numbered from 1 to 4, and we want to buy a contiguous area. Let $x_i$ be a binary variable denoting if we buy lot $i$.\footnote{This problem is inspired by a modeling problem in the author's work~\cite{template}.}
\begin{enumerate}[(a)]
\item If we should always buy lot 1 in case we were to buy any lots, how can we model it?
\item In case we do not have to necessarily buy lot 1 if we buy something, two of the three sets of constraints below---(i), (ii), and (iii)---would correctly model the problem, and one of those strictly dominates the other. Which one is which? And why?
\begin{enumerate}[(i)]
\item $x_1 + x_3 - x_2 \leq 1$ \\ $x_2 + x_4 - x_3 \leq 1$
\item $x_1 + x_3 - x_2 \leq 1$ \\ $x_2 + x_4 - x_3 \leq 1$ \\ $x_1 + x_4 - x_2 \leq 1$ \\ $x_1 + x_4 - x_3 \leq 1$
\item $x_1 + x_3 - x_2 \leq 1$ \\ $x_2 + x_4 - x_3 \leq 1$ \\ $x_1 + x_4 - x_2 -x_3 \leq 1$ 
\end{enumerate}
\end{enumerate}
\end{exercise}

\section{Using Solver Callbacks}\label{sec:callback}

Going back to the TSP model in Section~\ref{sec:tsp}, 
the main challenge that we face as the problem grows in size 
is that the total number of subtour elimination constraints quickly becomes prohibitive. 
However, not all of those constraints are really necessary. 
For example, 
we may get away with ignoring the most expensive subtours, 
such as the subtour involving only New York and Los Angeles, 
in case a solution with that subtour is worse than at least one feasible solution. 
More generally, 
what matters in such cases is not that the constraints perfectly describe the feasible set,
but instead that the constraints that are relevant to finding an optimal solution are part of the model.  
One way to do that is by adding the constraints only as needed.

If we start with the TSP model having only the first two constraints, 
which actually corresponds to listing Two-Way Assignment Model, 
we can modify the code to add subtour elimination constraints 
only after we conclude that each of those constraints is indeed necessary. 
These fashionably late constraints are called \emph{lazy constraints}.  
They are added with \emph{callbacks}, 
which are functions called at relevant moments by the solver to allow us to intervene on what happens next when the problem is being solved. 
Note that this is more efficient than solving the model again after adding a missed constraint: 
It is as if we were pretending that the lazy constraint was always there, 
but we simply chose not to take it into account earlier. 

To avoid excessive repetition of code from previous listings, 
we only describe the additions and how to adapt the listing Two-Way Assignment Model. 
In the first listing below, Using Lazy Constraints, 
we have a replacement for the last two lines. 
As can be observed in those lines, 
the name of the model object has been changed from \texttt{twap\_model} to \texttt{tsp\_callback\_model}. 
Given the use of the \texttt{time} library in this listing, 
it is also necessary to import that library as shown in Line 2 of TSP Model.  \\

\begin{Verbatim}[numbers=left,xleftmargin=5mm,frame=single,label=Using Lazy Constraints,fontsize=\small]
start_time = time.time()
tsp_callback_model.params.LazyConstraints = 1 
tsp_callback_model.optimize(callback_function)
end_time = time.time()
visualize_solution(x, coordinates, distance)
print("Duration:", end_time - start_time)
\end{Verbatim}
Besides the runtime reporting, which is not entirely new, 
the difference is in Lines 2 and 3. 
In Line 2 we change a parameter of the solver 
for allowing new constraints to be added to the model while it is being solved. 
In Line 3 we call the solver with \emph{.optimize}, 
but this time we provide the callback function as an argument. 

The callback function is meant to be used when the solver finds a supposedly feasible solution. 
If that solution has a subtour, our callback function generates the corresponding subtour elimination constraint. 
We present it in listing TSP Callback Function, which can be inserted after the first two lines of the adaptation of Two-Way Assignment Model that we are discussing. \\
 
\begin{Verbatim}[numbers=left,xleftmargin=5mm,frame=single,label=TSP Callback Function,fontsize=\small]
def callback_function(model, where):
    global x
    global cities
    global distance
    if where == gb.GRB.Callback.MIPSOL:
        current_x = model.cbGetSolution(x)
        next = cities[0]
        tour = []
        while True:
            tour.append(next)
            for (i,j) in distance:
                if i == next and current_x[i,j] > 0.99:
                    next = j
                    break
            if next == cities[0]:
                break
        if len(tour) < len(cities):
            print(tour)
            model.cbLazy(
                gb.quicksum( x[i,j]
                     for i in tour
                     for j in cities if j not in tour
                   ) >= 1
            )\end{Verbatim}
Callback functions should be specified in \texttt{gurobipy} with two arguments: 
(1) \texttt{model}, which allows us to interact with the model currently being solved; and 
(2) \texttt{where}, which allows us to find out what exactly is happening at the point that the callback is running. 
We compare \texttt{where} with \texttt{gb.GRB.Callback.MIPSOL} 
in Line 5 to check if this call was triggered by a feasible solution being found or not. 
This is a little different from what other solvers do. 
For example, 
we can use different callback functions in the CPLEX solver to handle different events~\cite{cplex}.  

Note that we use methods with slightly different names to interact with the model while the model is being solved, 
such as \texttt{.cbGetSolution} in Line 6 for extracting the feasible solution just found 
and \texttt{.cbLazy} in Lines 19 to 24 for adding a new lazy constraint. 
Between Lines 7 and 16, 
the listing describes how to start from the first city in the list of cities, 
\texttt{city[0]}, 
and find out to which city we go next with the solution, 
and then repeat that with the next city until we reach the first city again. 
As a numerical precaution in Line 12, 
we check if the binary variables are greater than 0.99 instead of equal to 1.
From Line 17, 
we generate the lazy constraint only if the tour does not include all cities. 
Otherwise, there is no subtour to be removed. 

Note that we can also generate lazy constraints before finding an integer solution  
by figuring out if they can remove an optimal fractional solution obtained at each branch-and-bound node. 
We return to that idea in the last exercise of this section.

\begin{exercise}
Build a table with the number of cities $n$ and the runtime of the TSP models with and without lazy constraints. 
Run the models to find the number of seconds for $n = 10, 11,$ and so on. Stop when either model takes too long.
\end{exercise}

\begin{svgraybox}
To play safe when running models that may take too long, 
you can use the instruction \texttt{m.Params.TimeLimit = 600} to stop the solver after 10 minutes (600 seconds), or another time of choice, 
where \texttt{m} is the name of the mathematical model object in Python. 
If you do that, 
you will need to look carefully at the log to find out if you got a solution and if that solution was proven to be optimal.
\end{svgraybox}

\begin{exercise}
The callback function presented removes the subtour involving the first city of the instance. 
Create new versions of the callback function in which (i) all subtours are removed; (ii) only the smallest subtour is removed; and (iii) only the largest subtour is removed. 
Build a table to compare the runtimes of the three new variants with the original callback function as $n$ grows.
\end{exercise}

\begin{exercise}
Using the TSP model with all subtour elimination constraints, 
another possible use of the callback is to remove every feasible solution when it is found (using the canonical cuts on the unit hypercube from Section~\ref{sec:tricks}) 
to count the total number of feasible solutions of the problem. Implement this callback function and test it with sizes starting at $n=3$. 
Then try it with another binary model. 
Note that you will need a global variable to be incremented inside the callback, 
and that the solver will conclude that no feasible solution exists when you are done. 
\end{exercise}

\begin{exercise}
Can you create a new callback function for removing fractional solutions using subtour elimination constraints? 
Instead of \texttt{MIPSOL}, you will need to check if the argument \texttt{where} matches \texttt{MIPNODE}. 
In order to retrieve the solution of LP relaxation from the branch-and-bound node, 
you need to use the method \texttt{.cbGetNodeRel}. 
Please note that the values can be anywhere from 0 to 1. 
Hence, your goal is not to extract a subtour directly from the solution, 
but instead to determine if there is a set $S$ such that the sum of all the variables in all arcs departing $S$ is smaller than 1. 
If you find more than one set $S$, 
you can test different strategies: (i) adding the constraint for removing subtours with the set leading to the smallest sum smaller than 1; (ii) adding the constraint for removing subtours with the set leading to the largest sum smaller than 1; or (iii) adding the constraint for removing subtours with all the sets leading to all the sums that are smaller than 1.
\end{exercise}

\begin{exercise}\label{ex:scatter}
Use your favorite statistical visualization tool for creating a scatter plot to compare the runtimes for solving a same TSP instance with TSP Model and with the model using callbacks that we described in this section. 
In such a scatter plot, 
every point represents a different instance. 
The value in $x$ is the time to solve using TSP Model. 
The value in $y$ is the time to solve using callbacks. 

For every point below the identity line ($x = y$), 
we know that $x > y$ and thus TSP Model takes longer to solve the instance. 
Conversely for every point above the identity line, 
we know that $x < y$ and thus TSP model takes less time to solve the instance. 
One way to compare two implementations computationally is by evaluating which one solves more instances faster, 
such as by counting points below and above the identity line. 
Hence, drawing the identity line in the chart might be helpful. 
\end{exercise}

\section{Designing Your Own (Heuristic) Algorithms}\label{sec:heuristic}

No matter how strong your model is, 
no matter how smart your tricks are, 
at some point you will not be able to push the solver any further. 
That is when the dream of finding a timely and provably optimal solution fades away. 
But fear not: 
we can still find pretty good solutions at that point. 
In fact, 
that happens very often with large scale applications in areas such as logistics and production planning. 

What we do in these cases in design our own algorithm. 
However, 
we stop thinking about designing exhaustive search algorithms   
and think instead of simple strategies that can quickly produce a solution that is reasonably good. 
If we can find ways of generating many reasonably good solutions that are likely different from each other, 
or if we can find ways of doing small improvements to a solution once it is found, 
even better. 
We call these ``imperfect'' algorithms \emph{heuristics}. 
They can be either \emph{constructive heuristics} if they create a solution from scratch, 
or \emph{local search heuristics} if they start from an existing solution and produce a new solution.

One of the simplest constructive heuristics is the \emph{greedy heuristic}. 
The idea is to build a solution by making decisions based on their immediate impact. 
For example, 
we start from some city in the TSP, 
find the closest city not visited yet, 
add the arc between the two to the solution, 
and then repeat from the next city until we are back to the first city. 
This is cheaper because we never change our choices based on what happens next, 
at which point we may find out that a local decision led to a much worse solution in the end. 
This strategy is in direct contrast with the exhaustive reasoning behind the branch-and-bound algorithm, 
which would have gone like this: 
try going to each one of the cities not visited yet as a separate subproblem, 
provided that we still expect to find a better solution out of it, 
and then compare all of their solutions. 
The listing below, Greedy TSP Heuristic, 
shows the cheaper alternative: \\

\begin{Verbatim}[numbers=left,xleftmargin=5mm,frame=single,label=Greedy TSP Heuristic,fontsize=\small]
cities, coordinates, distance = read_problem(
    "1000-cities.csv", 1000)

first_city = cities[0]
solution = [ first_city ] 
cities_left = [ c for c in cities if c != first_city ]
while len(cities_left) > 0: 
    current_city = solution[-1] 
    cities_left.sort(key=lambda c : distance[current_city,c]) 
    next_city = cities_left[0] 
    solution.append(next_city) 
    cities_left.remove(next_city) 
solution_value = total_distance(distance, solution)
    
print(solution_value, solution) 
show(coordinates, solution) 
\end{Verbatim}
Note that this listing requires functions from both listings with code in the appendix. 

When running this code for the whole 1000 cities, 
you may feel a deep satisfaction in how fast this turns out to be. 
But, of course, everything comes with a price: 
not every choice for this tour will seem to be a reasonable one afterwards. 
Back in the 1970s and 1980s, when finding optimal solutions was very difficult even for reasonably sized problems, 
many interesting ideas came up about what to do instead.
In fact, these ideas remain alive and new ones are still being proposed because we can always find a larger or more complicated problem to solve.\footnote{ 
We discuss alternative approaches with historical references from when they were proposed as well as more recent references on 
their application. You may need them one day. Besides, this section is more generours in citations than the previous ones because we cover other approaches here. 
For more references on the topics covered before in the text, see the next (and last) section.}

\subsection{Classic Heuristic Strategies}

There are two obvious ways to avoid making bad choices. 
First, we can try generating lots of solutions to use the best one that we can find through \emph{constructive heuristics}. 
A sensible way that does not give up on making reasonable decisions at each step 
are semi-gredy heuristics~\cite{hart1987semigreedy,feo1989setcovering}, 
in which we randomly pick among the top-ranked options. 
For example, 
rather than always picking the closest unvisited city to move next in the TSP, 
we pick a random city among the three closest ones. 

\begin{figure}[h]
\centering
\includegraphics[width=\textwidth]{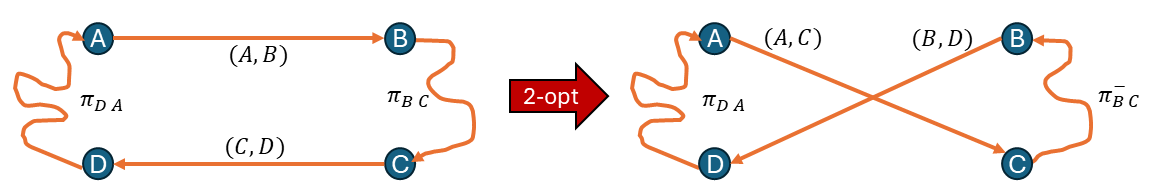}
\caption{Possible tour modification from applying 2-opt to obrain a different tour by replacing arcs $(A,B)$ and $(C,D)$ with $(A,C)$ and $(B,D)$, in addition to reversing the path $\pi_{B C}$ as $\pi^-_{B C}$.}
\label{fig:2Opt}       
\end{figure}

Second, 
as just mentioned above, 
we can simply try improving a solution based on what is obviously bad about it. 
This is where we think about a \emph{local search heuristic}. 
For example, 
suppose that we have a TSP solution in which we go directly from Los Angeles to New York without stopping by other cities along the way:
very likely, we will have a longer tour overall because we need to visit those other cities anyway. 
More generally, 
imagine that you have a very long TSP tour 
in which going directly from city $A$ to city $B$ 
seems to be a bad idea, 
and let $C$ and $D$ be another pair of consecutive cities in this tour. 
We would like to try connecting $A$ to one of those other cities, say $C$, and then $B$ is connected to $D$. 
If we represent the tour by the sequence of cities in it, 
say $A~B~\pi_{B C}~C~D~\pi_{D A}~A$, 
where $\pi_{i j}$ is a path with the cities immediately after $i$ and before $j$ in the tour, 
then with those changes we can get another tour $A~C~\pi^{-}_{B C}~B~D~\pi_{D A}~A$, 
in which $\pi^{-}_{B C}$ reverses the order in which we visit the cities that were immediately before $B$ and after $C$. 
If the cost of the new arcs used is lower, i.e., $d_{A C} + d_{B D} < d_{A B} + d_{C D}$, 
then this might be an improvement. 
However, 
if the direction matters, as in the assymetric TSP, 
then we also need to evaluate the cost of paths $\pi_{B C}$ and $\pi^{-}_{B C}$.
Figure~\ref{fig:2Opt} illustrates the swap of arcs and path reversal implied by the change. 
The 2-opt algorithm consists of doing such changes for as long as they improve the objective function~\cite{flood1956tsp,croes1958inversion}. 
The extension based on swapping three arcs is known as the 
3-opt algorithm~\cite{lin1965tsp}. 

\subsection{Metaheuristics}\label{sec:meta}

As no simple heuristic can be a silver bullet to solving NP-hard problems like the TSP, 
we are likely to converge to solutions that are not necessarily optimal. 
These are known as \emph{local optima}---as opposed to \emph{global optima}.\footnote{Optima is the plural of optimum.} 
To avoid stopping too far away from an optimal solution, 
more sophisticated and general-purpose techniques have been studied to explore the solution space. 
These techniques, which are then adapted to different optimization problems, 
are called \emph{metaheuristics}.

On the constructive side, 
we mention one cost-effective technique. 
GRASP (Greedy Randomized Adaptive Search Procedures)~\cite{feo1995grasp,Resende2010} 
combines the generation of semi-greedy solutions with local search on each generated solution.  
The classic application of GRASP is the \emph{set covering problem}. 
In this problem, we need to select sets numbered from $1$ to $m$ such that we have selected at least one set containing each item numbered from $1$ and $n$, with each set having a subset of the items from $1$ to $n$. 
A greedy way to solve this problem would be to pick the set with the most items not \emph{covered}\footnote{In set covering notation, an item is covered if we have selected a set that contains it.} yet, 
and then proceed by selecting other sets until all items are covered. 
We can randomize this greedy strategy by selecting one among the sets with more items not covered yet. 
By doing that, we can now repeat this construction multiple times to obtain different solutions for the problem. 
For each solution that we obtain, 
we may check if we can remove some of the sets selected if there are other selected sets covering all of their elements.

On the local search side, 
there has been a wide variety of approaches, 
some of which covered in detail in~\cite{local_search1,local_search2}. 
We mention three of the most popular here. 
Simulated Annealing~(SA) uses a decreasing tolerance on the change of the objective function, 
which at first allows the local search to move to much worse solutions, 
and then slowly comes to a point at which only strictly better solutions are acceptable. 
This is inspired on the physical process of annealing, which is applied by slowly cooling glass and metal in order to make the material more malleable~\cite{sa_kirkpatrick,sa_cerny,sa_now}. 
While also tolerating moving to a worse solution, 
Tabu Search~(TS) keeps a short-term memory of recently visited solutions or recently performed operations to avoid repeating them~\cite{ts1,ts2,ts3}. 
Another way to circumvent getting stuck into local optima is to work with many solutions at once. 
One of the most popular approachs is Genetic Algorithms~(GAs), 
in which we model the solutions of a problem as if they were a DNA sequence, 
such as the sequence of cities visited in TSP tours. 
Hence, different solutions represent different individuals in a population that keeps evolving: 
New solutions are obtained when individuals mutate or produce offspring, 
and solutions are discarded by subjecting inviduals to survival criteria that relate to the objective function value~\cite{ga1,ga2}. 

It is important to observe that the literature on metaheuristics is both very creative and very extensive. 
However, 
there is a growing and vocal concern among scholars that some of this work is more centered on the metaphors 
than on actually making progress in solving larger and harder optimization problems~\cite{metaphor1,metaphor2}.

\subsection{Quantum Computing and Machine Learning}

More recently, 
the ideas from SA have been applied for designing specialized quantum computers, 
which can be used to model discrete optimization problems as quantum systems in which an optimal solution corresponds to the state of lowest energy of the system~\cite{Farhi2000,Farhi2001,Lucas2014}. 
Note that we may not be able to solve very large problems with such machines yet, 
and that their solutions may be within a certain gap from the optimal value~\cite{Dash}. 
This is one among many quantum computing paradigms currently in study, 
which may one day change how we approach large scale optimization problems~\cite{nature,nasa}. 

There is also a growing use of machine learning as a means to improve how we solve optimization problems. 
This has been a very rich area of work in recent years, which has been covered by surveys such as~\cite{AdaptSurvey,CombOptTour,OptimizationSurvey}. 
Most of it would not be of immediate concern for someone starting to learn and use optimization, 
except perhaps for \emph{optimization proxies}. 
An optimization proxy is a machine learning model, which is trained on a set of instances and their solutions, 
to quickly predict a good solution for a new instance without using an optimization algorithm for doing so. 
This may lead to solutions that are infeasible and require some adjustments, 
but they have been shown to be effective in solving certain problems faster than an optimization solver would be able to handle them~\cite{proxies1,proxies2,proxies3}. 

\begin{exercise}
If you look carefully at the result of running the listing Greedy TSP Heuristic, 
you will notice some decisions that are obviously bad in retrospect---especially towards the end. 
Which decisions are those? Why are they bad? 
\end{exercise}

\begin{exercise}
To produce different solutions with the greedy heuristic, 
we can change the city at which we start constructing the tour. 
\begin{enumerate}[(a)]
\item Why would this produce a different tour? 
\item Make a new version of the listing in which you try using each city as the starting point and then pick the best tour produced. 
You may consider limiting the number of cities $n$ to 100 to avoid excessive runtime. 
\end{enumerate}
\end{exercise}

\begin{exercise}
Implement the 2-opt algorithm as a loop that keeps looking for pairs of arcs that could be replaced to reduce the cost of the tour. 
Add this algorithm as a follow-up step to the greedy algorithm. 
By how much does it improve the total distance traveled in a tour?
\end{exercise}

\begin{exercise}
We have seen that the 2-opt algorithm requires reversing a path in the tour, 
which may entail recalculating the distance traveled along the reversed path if the TSP is assymetric. 
By thinking of 3-opt as a modification in which three arcs are removed and their ends are reconnected, 
answer the following questions: 
\begin{enumerate}[(a)]
\item It is possible to use 3-opt without having to reverse any paths?
\item How many tours can be produced from a given tour using 3-opt?
\end{enumerate}
\end{exercise}

\chproblem{
Implement a random generator of instances for the set covering problem and an ILP model to solve them. 
Then implement the GRASP-based heuristic described in Section~\ref{sec:meta} to solve those instances. 
Compare the runtimes and the value of the best solutions found by solving the model vs. running the heuristic multiple times on instances with varying numbers of sets and items. 

Consider the following questions: 
Can you find an optimal solution with GRASP? How many runs does it take to get there? Is that faster than solving the ILP model?

See the instructions in other challenge problems and their footnotes for guidance on how to conduct similiar experiments and implement randomization in Python. 
You may also consider using scatter plots as in Exercise~\ref{ex:scatter}, provided that the comparisons presented in such plots are reasonable.\footnote{
This is a situation in which we risk ``comparing apples and oranges'' if we are not careful, in the sense that the results are not meaningful. 
Comparing the runtime of solving an ILP model to optimality with the time to generate a solution that is not optimal would not make sense. However, you may consider comparing the former with the time by which the first optimal solution is generated with an heuristic (even if we cannot tell if a solution is optimal by using the heuristic alone, in this case we would know that from having solved the ILP model and that helps us evaluate how good the heuristic is). From a computational perspective, that would tell us if putting the same amount of computational effort would possibly produce an optimal solution with an heuristic.}
}

\section{Where to Go Next}

We have seen that solving optimization problems can be as simple as just modeling them, 
but in many cases that will not be enough: 
we may need to (i) analyze multiple forms of modeling; 
(ii) consider ways of intervening in how the solver handles the model; 
or even (iii) give up on optimality and write our own heuristics from scratch.

Ultimately, 
most optimization problems worth solving are problems that someone already has to deal with on their own. 
Sometimes, finding a sufficiently good solution in a reasonable amount time would already be a great win 
in comparison to how the problem would have been solved otherwise. 
For that reason, 
trying to replicate in code what a decision maker currently does is a great first step, 
since it gives us a benchmark to evaluate our own work later. 
Besides, 
improvements to the algorithms within solvers and to the hardware in which we use them 
keep pushing the lines of what problems---and sizes---would require the sophistication of (i), (ii), or (iii)~\cite{bixby2012}.

For the students interested in diving further on discrete optimization, 
I would recommend a combination of theoretical grounding with hands-on experience. 
Following a broader introduction to discrete optimization with a book such as~\cite{co_book}, 
there are many books focused on the broader area of Integer Programming~(IP)~\cite{bertsimas,ccz,schrijver,wolsey,nemhauser}. 
There is not one better or more complete IP book:  
Different authors focus on the parts that are closer to their interests and that are more relevant to how they see the subject. 
There are also relevant books and surveys in mathematical modeling for optimization~\cite{aimms,mosek,model}; 
strengthening models by adding valid inequalities, or cuts, and the underlying polyhedral theory supporting the generation of such cuts~\cite{balas,inequalities,tuncel,ziegler}; 
alternative modeling techniques such as constraint programming~\cite{cbs,dechter,cp_handbook} and decision diagrams~\cite{bdd,bdd_survey}; 
optimization under uncertainty~\cite{stochastic}; 
integrating different optimization methods~\cite{integrated}; 
embedding neural networks in optimization models~\cite{huchette2023survey};  
optimization in market design~\cite{bichler,vohra};
the traveling salesperson problem~\cite{tsp_book}\footnote{For a more accessible introduction to the TSP, see~\cite{tsp_pop}.}; 
and using linear optimization for art~\cite{opt_art}. 

Based on the exercises and challenge problems that we discussed, 
here are some specific ideas for research projects that you could try on your own:

\resproject{\textbf{Heuristics:}
Design an heuristic for solving an optimization problem that interests you, 
and then find out how large the problem should be for your heuristic to be a viable option in comparison to writing an ILP model directly. 
For example, the problem may get so large that it may take less time to find a good --- but possibly not optimal --- solution with the heuristic. 

You may find it easier to start by using one of the problems and one of the heuristics that we have discussed first, 
and then gradually move on to experimenting with a different problem or heuristic afterwards.   
}

\resproject{\textbf{Branch-and-Bound:} 
Implement your own branch-and-bound algorithm for solving discrete optimization problems by relying only on the LP solver from Gurobi, 
and then compare it with using the ILP solver from Gurobi. 
Very likely, your code will be a lot slower. 
However, you will have to think about a lot of details in your implementation, 
such as which decision variable to use when branching and which node to explore next in the branch-and-bound tree. 
By evaluating how those different decisions affect the runtime of different problems that you try to solve, 
you will grow your intuition for what makes a problem or its instances more difficult to solve. 

You may consider comparing the strategy of selecting the most fractional variable for branching with choosing a fractional variable at random, 
and then think about how to influence your choice based on the coefficient of each decision variable in the objective function. 
Eventually, you may find it relevant to look up what a reduced cost is in a book about linear programming, 
and then try to obtain the reduced costs after solving each LP relaxation in your code. 
}

\resproject{\textbf{Benchmark Problems:} 
There are many benchmarks of optimization problems and instances that are used for research. 
Some of those have problems for which we do not know the optimal value, 
or we know it but it takes a long time to solve those problems. 
Consider developing an algorithm that can find a good solution for one of such problems or instances. 
As a starting point on benchmarks and a tool to search for instances, see \cite{miplibing}.
}

Here are more general ideas on how to get hands-on experience and join projects:

\resproject{\textbf{Applying for Research Assistant Internships:}

Some universities offer summer research programs for their own students and also for students from other institutions. 
For local opportunities, 
look for advertisements of information sessions and campus jobs; 
attend events about research opportunites and in which students present their research; 
talk with professors who are teaching your courses; 
and reach out to other professors who work in what interests you. 
For mathematical optimization, 
there might be related opportunities in a variety of departments, 
such as Business Analytics, Computer Science, Engineering (Chemical, Computer, Civil, Electrical, Industrial, Mechanical), Mathematics (Applied, Pure), Information Systems, Management Science, Operations Research, Physics, and so on---it is possible that the right research experience for you might not be in the department where you study. 
In fact, some students switch departments in graduate school, and then switch back or switch again when they become professors.

For external opportunites, 
look for programs such as Research Experiences for Undergraduates (REU), 
which is supported by the National Science Fountation (NSF). 
That particular program supports universities in running special programs on a specific subject. 
You may also find some opportunities working for companies during the summer, 
but that may vary a lot by region. }

\resproject{\textbf{Tackling a Campus or Workplace Problem:}

You may find something interesting to optimize around you, 
such as common problems in academic campuses and workplaces: 
\begin{itemize}
\item Scheduling midterms and final exams outside lecture times is a common problem in many universities. 
Each university also has their very specific constraints, 
which may prevent using a tool developed elsewhere.
The authors has been involved in one of such projects at Bucknell University~\cite{febu}.
\item Scheduling work shifts at a library, restaurant, or help desk. 
Every worker has days and times in which they would prefer to work, or that they cannot work due to other commitments.
\item Carpooling assigments for commuting to work. 
Each person should be assigned to a group. Every group should have at least one driver. We want to minimize the extra time required from the driver.
\end{itemize}
}

\resproject{\textbf{Participating in an Optimization Challenge:}

Some academic conferences and websites organize competitions every year, 
at least in recent years as of the time of writing this. 
The references below are among the most recent occurrences of some of these: 
the MOPTA Conference organizes an applied  modeling competition~\cite{mopta};  
the MIP Workshop organizes a computational competition focused on improving specific aspects of optimization solvers~\cite{mip}; 
and some of the Kaggle competitions entail using optimization, such as the recurring Santa Challenge~\cite{kaggle}.}

\resproject{\textbf{Replicating an Academic Paper:}

Find a paper about a problem that interests you. Based on the description of the problem, try coming up with your own model. Then try implementing the model described in the paper, if that model is not too complicated. If there is data, try running either model with it. If there is no data, try to come up with your own data. You may find some interesting applications in academic journals such as the INFORMS Journal on Applied Analytics (formerly known as Interfaces), the European Journal of Operational Research, and International Transactions in Operational Research; 
and conferences such as CPAIOR.

There are also some journals focused on undergraduate mathematics, 
which may provide more accessible examples.  
The SIAM Undergraduate Research Online (SIURO) Journal~\cite{siuro} publishes articles authored by undergraduate students, 
including optimization models for planning the visit to a theme park~\cite{disney} and the expansion of mental health facilities~\cite{mentalHealth}, as well as metaheuristics like genetic algorithms applied to assignment problems~\cite{qapGenetic} and to flip a beefsteak~\cite{beefsteak}. 
The Undergraduate Mathematics and Its Applications (UMAP) Journal~\cite{umap} is a modeling journal that provides access to their latest issue with a free membership, 
with a recent issue featuring a student project on optimizing the assignment of first-year seminars~\cite{seminar}. 
This journal is part of the Consortium for Mathematics and Its Applications (COMAP), 
which also publishes modeling modules related to optimization, 
such as for airline scheduling~\cite{airline} and for prison guard scheduling~\cite{prison}. 
The Involve Journal~\cite{involve} is also focused on the participation of students in mathematical research and all issues can be freely downloaded.  
One interesting paper presents a new optimization algorithm developed in a course project~\cite{mst}. 
The journal also has articles on career advice, 
including a piece on nonacademic careers that has a paragraph about operations research~\cite{nonacademic} 
and another piece on where to find ideas and data for sports analytics projects~\cite{sports}.
}

In case you succeed in developing a research project as an undergraduate student and write a report about it, 
consider submitting it to one of the journals above and to competitions like the INFORMS Undergraduate Operations Research Prize~\cite{urop}.

You may also find useful to attend academic seminars, 
in particular if they are not overly technical by the description. 
Similar to the case of research opportunities, 
you may find relevant seminars in many different departments on campus. 
In addition,
there are some series of online seminars that may capture your interest, 
such as Discrete Optimization Talks~\cite{dot}, 
the EURO Online Seminar Series on Operational Research and Machine Learning~\cite{euro}, 
the Gurobi Academic Webinar~\cite{gurobi_aw}, and 
the Robust Optimization Webinar~\cite{row}. 
Another useful resource are the recordings of the workshop Computational Optimization at Work, 
which covers in great depth how optimization solvers work and how to use them in applications~\cite{cowork}. 
Finally, 
you may find it interesting to listen to the trajectory of established scholars and professionals in the area of optimization 
by listening to the ``Subject to'' podcast~\cite{subjectTo}. 

By reading this far, I hope you have seen how fascinating optimization is. 
I wish you the best in your optimization journey. 
Actually, I wish you \emph{an} optimal journey!


\section*{Appendix}
\addcontentsline{toc}{section}{Appendix}

\section*{Supporting Functions for Listings in Sections~\ref{sec:tsp}, \ref{sec:callback}, and \ref{sec:heuristic}}

\begin{Verbatim}[numbers=left,xleftmargin=5mm,frame=single,label=Appendix Code 1,fontsize=\small]
import math
import matplotlib.pyplot as plt
import networkx as nx
import warnings

from itertools import chain, combinations

def __distance_between_cities(coordinates: dict, A, B):
    (latA,lonA) = coordinates[A]
    (latB,lonB) = coordinates[B]
    return 62.36*math.sqrt(
	(latA-latB)*(latA-latB) + (lonA-lonB)*(lonA-lonB))

def read_problem(file_name: str, size: int):
    f = open(file_name) 
    cities = []
    coordinates = {}
    distances = {}
    for i in range(size): # Loop for reading each city
        line = f.readline()
        tokens = line.split(",")
        city_name = tokens[1]
        city_coord = (float(tokens[2]), float(tokens[3]))
        cities.append(city_name)
        coordinates[city_name] = city_coord
    for i in cities:
        for j in cities:
            distances[(i,j)] = 
			__distance_between_cities(
				coordinates, i, j)
    return cities, coordinates, distances

def visualize_solution(x, coordinates, distances):
    warnings.filterwarnings("ignore", category=UserWarning)
    G = nx.DiGraph()
    for city in coordinates: 
        G.add_node(city)
    for (i,j) in distances:
        if x[i,j].X > 0.99:
            G.add_edge(i,j)
            print(i, "->", j)
    nx.draw(G, pos=coordinates, with_labels=True, 
		font_weight='bold', edge_color="red", width=3)
    plt.margins(x=0.1)
    plt.show()

def subtours(cities):
    min_size = 2
    max_size = len(cities) // 2
    return chain.from_iterable(
        combinations(cities, r) 
        for r in range(min_size, max_size+1))
\end{Verbatim}

\section*{Supporting Data for Models}

The listing below contains the first 20 lines of the CSV file used for the experiments with TSP models, 
which is provided to make this piece self-contained. 
For the entire dataset, 
visit \url{https://github.com/thserra/discreet}.
This dataset was derived from \cite{opendatasoft}, 
which as of the writing is no longer available online. 

\begin{Verbatim}[numbers=left,xleftmargin=5mm,frame=single,label=Appendix TSP Data,fontsize=\small]
1,New York-NY,-74.0059413,40.7127837
2,Los Angeles-CA,-118.2436849,34.0522342
3,Chicago-IL,-87.6297982,41.8781136
4,Houston-TX,-95.3698028,29.7604267
5,Philadelphia-PA,-75.1652215,39.9525839
6,Phoenix-AZ,-112.0740373,33.4483771
7,San Antonio-TX,-98.4936282,29.4241219
8,San Diego-CA,-117.1610838,32.715738
9,Dallas-TX,-96.7969879,32.7766642
10,San Jose-CA,-121.8863286,37.3382082
11,Austin-TX,-97.7430608,30.267153
12,Indianapolis-IN,-86.158068,39.768403
13,Jacksonville-FL,-81.655651,30.3321838
14,San Francisco-CA,-122.4194155,37.7749295
15,Columbus-OH,-82.9987942,39.9611755
16,Charlotte-NC,-80.8431267,35.2270869
17,Fort Worth-TX,-97.3307658,32.7554883
18,Detroit-MI,-83.0457538,42.331427
19,El Paso-TX,-106.4424559,31.7775757
20,Memphis-TN,-90.0489801,35.1495343
\end{Verbatim}


\section*{Additional Supporting Functions for Listing in Section~\ref{sec:heuristic}}

The following listing should be used in addition to the listing Appendix Code 1, 
in which all required libraries are imported.

\begin{Verbatim}[numbers=left,xleftmargin=5mm,frame=single,label=Appendix Code 2,fontsize=\small]
def total_distance(distance, solution): 
    t_distance = 0
    N = len(solution)
    for i in range(N-1): 
        t_distance += distance[solution[i],solution[i+1]]
    t_distance += distance[solution[N-1],solution[0]] 
    return t_distance

def show(coordinates, solution): 
    G = nx.DiGraph()
    for city in cities: 
        G.add_node(city)
    N = len(solution)
    for i in range(N-1): 
        G.add_edge(solution[i],solution[i+1])
    G.add_edge(solution[N-1],solution[0]) 
    nx.draw(G, pos=coordinates, with_labels=True,
		font_weight='bold', edge_color="red", width=3)
    plt.margins(x=0.1)
    plt.show()
\end{Verbatim}

\bibliographystyle{splncs04}
\bibliography{references}

\end{document}